\documentclass[11pt,twoside,a4paper]{article}
\usepackage[english]{babel}
\usepackage{a4wide}

\usepackage{epsf,epsfig}
\usepackage{graphics,color}
\usepackage{amsmath}
\usepackage{amssymb}
\usepackage{cite}
\usepackage{verbatim}
\usepackage{float}
\usepackage{graphicx}
\usepackage{amsthm}
\usepackage{textcomp}
\usepackage{subfig}
\usepackage{framed}
\usepackage{url}
\usepackage{tikz}
\usepackage{enumerate}
\usepackage{setspace}

\newcommand{\curl}{\operatorname{curl}}
\newcommand{\dist}{\operatorname{dist}}

\newcommand{\qtext}[1]{\quad \text{#1} \quad}
\newcommand{\qand}{\qtext{and}}

\newcommand{\xto}[1]{\xrightarrow{#1}}

\newcommand{\e}{\varepsilon}
\newcommand{\N}{\mathbb N}
\newcommand{\R}{\mathbb R}

\newcommand{\Z}{\mathbb Z}

\newcommand{\bx}{\mathbf x}

\newcommand{\by}{\mathbf y}

\newcommand{\cC}{\mathcal C}

\newcommand{\cF}{\mathcal F}

\newtheorem{thm}{Theorem}[section]

\newtheorem{rem}[thm]{Remark}

\title{Expansions for the linear-elastic contribution to the self-interaction force of dislocation curves}
\author{Patrick van Meurs}

\begin{document}


\maketitle



\begin{abstract}
The self-interaction force of dislocation curves in metals depends on the local arrangement of the atoms and on the nonlocal interaction between dislocation curve segments. While these nonlocal segment-segment interactions can be accurately described by linear elasticity when the segments are further apart than the atomic scale of size $\e$, this model breaks down and blows up when the segments are $O(\e)$ apart. To separate the nonlocal interactions from the local contribution, various models depending on $\e$ have been constructed to account for the nonlocal term. However, there are no quantitative comparisons available between these models. This paper makes such comparisons possible by expanding the self-interaction force in these models in $\e$ beyond the $O(1)$-term. Our derivation of these expansions relies on asymptotic analysis. The practical use of these expansions is demonstrated by developing numerical schemes for them, and by -- for the first time -- bounding the corresponding discretization error.
\end{abstract}

\section{Introduction}

Dislocations in metals are curve-like defects in the atomic lattice. The emergent group behavior of many dislocations is the driving mechanism of plastic deformation of metals. One of the reasons that describing plastic deformation is an active field of research is that there is no consensus on a description for the self-interaction force of a dislocation curve. Such a description has been sought for more than half a century; we cite several papers in the remainder of the introduction. The aim of this paper is to contribute towards reaching consensus by developing a framework for computing expansions for several descriptions of the self-interaction force. 

\begin{figure}[ht]
\centering
\begin{tikzpicture}[scale=1.5, >= latex]
\def \rr {0.04} 

\draw[->] (-1.2,1.5) --++ (1,0) node[below] {$e_1$};
\draw[->] (-.7, 1) --++ (0,1) node[left] {$e_2$};
\draw[dotted] (-1,0) -- (1,0);

\draw [thick] plot [smooth cycle, tension=.7] coordinates {(0,0) (0,1) (1,1.5) (1,.5) (1.5, 0) (1, -0.732)};

\draw[blue] (.2,0) node[anchor = south west]{$\phi(x)$} arc (0:210:.2) node[anchor = north west, black] {$x$};
\draw[dotted, rotate = 300] (-1,0) -- (1.5,0);
\draw[->, red, rotate = 300] (0,0) -- (1,0) node[anchor = north east] {$\tau(x)$};
\draw[->, red, rotate = 210] (0,0) -- (1,0) node[left] {$n(x)$};
\draw (1, -0.732) node[anchor = north west] {$\Gamma$};
\draw[<-] (1.1,.5) -- (2, .8) node[right] {$\kappa > 0$};
\draw[<-] (1.1,1.6) to [out=45,in=180, looseness=1] (2, 2) node[right] {$\kappa < 0$};
\fill (0,0) circle (\rr);

\begin{scope}[shift={(5,0)},scale=1] 
\draw [thick] plot [smooth cycle, tension=.7] coordinates {(0,0) (0,1) (1,1.5) (1,.5) (1.5, 0) (1, -0.732)};
\fill[white] (0,0) circle (.5);
\draw [dotted] plot [smooth cycle, tension=.7] coordinates {(0,0) (0,1) (1,1.5) (1,.5) (1.5, 0) (1, -0.732)};


\draw (0,0) circle (.5);
\draw (0,-.5) node[anchor = north east] {$B(x,\e)$};
\draw (1, -0.732) node[anchor = north west] {$\Gamma_\e(x)$};  
\fill (0,0) circle (\rr);      
\end{scope}
\end{tikzpicture} \\
\caption{Left: sketch of $\Gamma$. Right: sketch of $\Gamma_\e(x)$.}
\label{fig:Gamma0}
\end{figure}
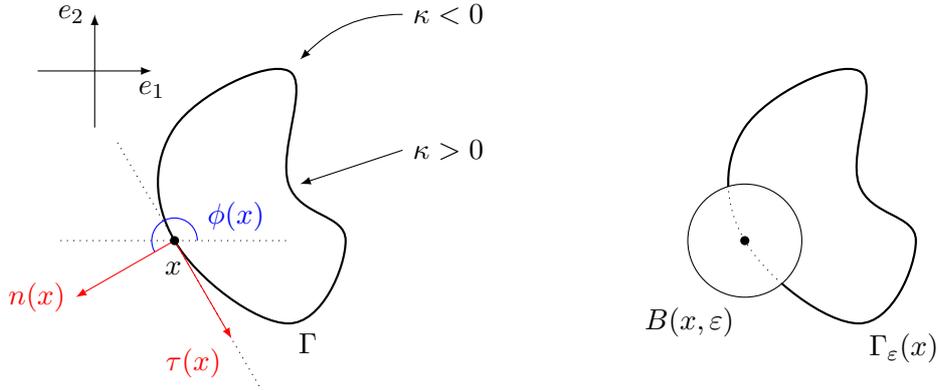

To describe the complexity of describing the self-interaction force of a dislocation curve, we first introduce the setting. To avoid long formulas and to keep the focus on the methodology, we consider the simple setting where $\R^2$ represents an isotropic elastic medium with shear modulus $\mu > 0$ and Poisson ratio $\nu \in (-1, \frac12)$. We consider a dislocation loop represented by a closed regular curve $\Gamma \subset \R^2$, and set $b = e_1$ as its Burgers vector, where $\{e_1, e_2\}$ is the standard basis in $\R^2$. For readers unfamiliar with dislocations and their Burgers vector, we refer to the textbooks \cite{HirthLothe82,HullBacon11}. For $x \in \Gamma$ we set
\begin{equation} \label{ktaun0}
  \kappa(x), \quad
  \tau(x)  = \begin{bmatrix}
  -\sin \phi(x) \\
  \cos \phi(x)
  \end{bmatrix}, \quad
  n(x) = \begin{bmatrix}
  \cos \phi(x) \\
  \sin \phi(x)
  \end{bmatrix} 
\end{equation} 
as respectively the curvature, tangent vector (counter-clockwise direction) and outward pointing normal vector of $\Gamma$ at $x$. Figure \ref{fig:Gamma0} illustrates the setting. Our sign convention is that $\kappa(x) < 0$ when $\Gamma$ is a circle. 

Next we describe the complexity of the self-interaction force on the dislocation $\Gamma$ at $x$. For practical applications such as dislocation dynamics, it is sufficient to focus on the (scalar-valued) normal component of this force. For brevity, in the remainder we simply refer to this normal component as the self-interaction force. The reader familiar with the self-interaction force will recognize formulas \eqref{F}--\eqref{Ge} below, and can skip the following explanatory text up to \eqref{Ge}. 

Instead of the self-interaction force on $\Gamma$ at $x$, let us first consider the easier expression for the force that a second dislocation loop $\tilde \Gamma$ with the same Burgers vector $b$ exerts on $\Gamma$ at $x$. We assume that $\tilde \Gamma$ is far enough away from $\Gamma$ such that the interaction force is accurately described by linear elasticity. Then, the force $\tilde F(x)$ exerted by $\tilde \Gamma$ on $\Gamma$ at $x$ is found by first computing the stress in the medium at $x$ caused by $\tilde \Gamma$ (see, e.g.\ \cite[(7.4)]{DeWit60}), and then by applying the Peach-Koehler formula \cite{PeachKoehler50} to convert this stress into a force. This yields
\begin{equation} \label{F}
  \tilde F(x)
  := \oint_{\tilde \Gamma} G(y-x) \cdot \tau (y) \, dy,
\end{equation}
where we have normalized the multiplicative constant to $1$, and
\begin{equation*} 
   G(z) := \left( \frac{z^T}{ |z|^3 } \begin{bmatrix}
  0 & -1 \\
  1-\nu & 0
  \end{bmatrix} \right)^T = \frac1{ |z|^3 } \begin{bmatrix}
  (1-\nu) z_2 \\ - z_1
  \end{bmatrix}.
\end{equation*}
Due to the singularity of $G(z)$ at $z = 0$, the function $\tilde F : \R^2 \to \R$ is singular on $\Gamma$. This stems from the fact that linear elasticity is a poor model for describing the stress field in the \emph{dislocation core}. The dislocation core is the tubular neighborhood with radius of about $5$ atoms ($\sim 1 $ nm) around a dislocation. In particular, this means that \eqref{F} cannot be used to describe the force that $\Gamma$ exerts on itself at any point $x \in \Gamma$.
 
However, in practice, the width of the dislocation core is much smaller than the diameter of $\Gamma$, which we take for the sake of explanation to be the typical value of $10 \mu$m. Then, scaling the diameter of the dislocation loop to $1$ dimensionless unit, the radius of the dislocation core is
\begin{equation} \label{e:10m4}
  \e \sim 10^{-4}.
\end{equation}
In the remainder we regard $\e$ as a generic small parameter which indicates the length scale where linear elasticity breaks down. We leave the choice of $\e$ as a modeling issue.

The easiest modeling choice for avoiding the breakdown of linear elasticity inside the dislocation core is to neglect the contribution of the dislocation core to the self-interaction force. To describe the resulting force, we set
\[
  \Gamma_\varepsilon(x) := \Gamma \setminus B(x, \varepsilon);
\]
see Figure \ref{fig:Gamma0} for a sketch. Then, \eqref{F} implies that
\begin{equation} \label{Fe}
  F_\e(x)
  := \int_{\Gamma_\varepsilon(x)} G(y-x) \cdot \tau (y) \, dy
\end{equation}
gives an acceptable approximation of the force exerted on $\Gamma$ at $x$ by $\Gamma_\e(x)$. This model can be made more accurate by coupling it to an atomistic model which accounts for the local contribution of $\Gamma \cap B(x, \e)$; see e.g.\ 
\cite{GehlenHirthHoaglandKanninen72,
HenagerHoagland05,
Clouet11}.

In this paper we focus on the contribution from linear elasticity to the self-interaction force on $\Gamma$. Since there is no clear splitting of the self-interaction force in terms of a linear-elastic part and microscopic contributions, there are several choices on how to define the linear elastic part. \eqref{Fe} is one such choice; it dates back at least to \cite{HirthLothe82}. Other choices are derived by first smearing out the singularity of the stress field within the dislocation core, and then by computing the force from this regularized stress field. This is a popular approach, because it includes a phenomenological model for the dislocation core, which then needs no additional treatment. The smearing out of the dislocation core has been done either with a (singular) convolution kernel (see \cite{CaiArsenlisWeinbergerBulatov06} for a short review), or with a phase-field approach (see \cite{KundinEmmerichZimmer11} for a review). In particular, we focus on the specific convolution kernel constructed in \cite{CaiArsenlisWeinbergerBulatov06}. The feature of this convolution kernel is that the convolution integral in the related expression for $G$ can be computed explicitly; see \cite[(33)]{CaiArsenlisWeinbergerBulatov06}\footnote{Precisely, we take the regularized stress field induced by a dislocation, and use it in the formula for the Peach-Koehler force applied to a (singular) dislocation.}. The resulting self-interaction force is given by  
\begin{align} \label{cFe} 
  \cF_\e(x)
  = \oint_\Gamma G_\e(y-x) \cdot \tau (y) \, dy,
\end{align}
where
\begin{align} \label{Ge}
  G_\e(z)^T := \frac{z^T}{ \sqrt{ |z|^2 + \varepsilon^2 }^3 } \begin{bmatrix}
  0 & -1 \\
  1-\nu & 0
  \end{bmatrix}
  + \frac{3 \varepsilon^2 z^T}{ 2 \sqrt{ |z|^2 + \varepsilon^2 }^5 } \begin{bmatrix}
  0 & 0 \\
  1 - \nu & 0
  \end{bmatrix}.
\end{align}
The difference with \eqref{Fe} is that the dependence on $\e$ in the integral in \eqref{cFe} has shifted from the integration domain to the integrand.

To summarize the above, $F_\e(x)$ and $\cF_\e(x)$ are two of the many choices for the (linear-elastic contribution of the) self-interaction force on $\Gamma$ at $x$.
While attempts have been made to compare different choices (see, e.g.\  \cite{LeSar04,CaiArsenlisWeinbergerBulatov06}), the comparisons mainly rely on unquantifiable notions and practical experiences when using them in numerical computations. To quantify the differences, one needs expansions of the self-interaction force in terms of $\e$. Constructing such expansions is the first of the two aims of this paper. The second aim is to use these expansions to construct accurate schemes by which they can be computed numerically. We pursue these two aims respectively in Sections \ref{s:i:ana} and \ref{s:i:num}.

\subsection{The expansions of $F_\e$ and $\cF_\e$}
\label{s:i:ana}

Expansions for forces such as $F_\e$ and $\cF_\e$ in terms of $\e$ have been constructed in 
\cite{GavazzaBarnett76, 
Lothe92,
ZhaoWangXiang2012} 
(each paper considers a different core regularization; none of which equals $F_\e$ or $\cF_\e$), and in a more general setting than in this paper (either in 3 dimensions or for an anisotropic medium). In all three works, the self-interaction force (reduced to our setting) expands as
\begin{equation} \label{F:expa:O1}  
  \kappa(x) \big(  1 + \nu - 3 \nu \cos^2 \phi(x) \big) \log \frac1\e
  + O(1)
\end{equation}
as $\e \to 0$. The leading order term of $O(|\log \e|)$ is consistent in all works, and also coincides with the term obtain from phase-field models. 

However, it is desired to specify the $O(1)$ term too. Indeed, for the typical value of $\e$ in \eqref{e:10m4}, the $O(1)$ term may not be significantly smaller than the $O(|\log \e|)$ term. In particular, for parts of $\Gamma$ which are approximately straight (i.e.\ where $|\kappa|$ is small), the prefactor of the $O(|\log \e|)$ term is small, whereas the nonlocal contribution to the interaction force (which contributes to the $O(1)$ term) does not decay (in general) in $\kappa$ for small $|\kappa|$. 

Yet, the treatment and expression of the $O(1)$ term is quite different in the three works cited above. In \cite[(156)]{Lothe92} this term is left unspecified. In \cite[(7.1)]{GavazzaBarnett76} only a part of this term is specified. In \cite[(44)]{ZhaoWangXiang2012} the self-interaction force of the curve segment $\Gamma \cap B(x, \e')$ at $x$ is expanded up to an error of $O(\e' + \e/(\e')^2)$, where $\e'$ is a phenomenological mesoscopic length scale which satisfies $\sqrt \e \ll \e' \ll 1$. The related $O(1)$ term is explicit, but local (i.e.\ it depends on $\Gamma$ only through $n(x)$ and $\kappa(x)$ at $x$). The contribution of $\Gamma_{\e'}$ on the self-interaction force is not expanded.

Hence, the results in the literature that go beyond an unspecified $O(1)$ term in \eqref{F:expa:O1} are very limited. Our first main contribution (Theorem \ref{t:ana}) extends these results by characterizing completely the $O(1)$ term in the expansions of $F_\e(x)$ and $\cF_\e(x)$, and by showing that the next order term is $O(\e)$. To state this result, we set 
\[
  \cC(x) := \left\{ \begin{aligned}
    &\{ t \tau(x) : t \in \R \}
    &&\text{if } \kappa(x) = 0 \\
    &\partial B \Big( x + \frac{ n(x) }{ \kappa(x) }, \frac1{|\kappa(x)|} \Big) 
    &&\text{otherwise}
  \end{aligned} \right.  
\]
as the tangent circle of $\Gamma$ at $x$. Similar to $\Gamma_\e(x)$, we set 
\[
 \cC_\e(x) := \cC(x) \setminus B(0, \e).
\]

\begin{thm}[Expansions of $F_\e(x)$ and $\cF_\e(x)$] \label{t:ana}
Let $\Gamma$ be a non-intersecting closed curve of class $C^3$ with finite length. Then, for all $x \in \Gamma$, the following limit exists
\begin{equation} \label{Psi:formal}
  \Psi(x) := \lim_{\e \to 0} \int_{\Gamma_\e(x) \cup \cC_\e(x)} G(y-x) \cdot \tau (y) \, dy
\end{equation}
and
\begin{align} \label{Fe:t}
  F_\e(x)
  &= \kappa(x) A_{\phi(x)} \log \frac1 {\varepsilon |\kappa(x)|} + \kappa(x) B_{\phi(x)} 
  + \Psi(x) + O (\varepsilon) \\\label{cFe:t}
  \cF_\e(x) &= F_\e(x) + \kappa(x) C_{\phi(x)} + O (\varepsilon)
\end{align}
as $\e \to 0$, where
\begin{subequations} \label{ABC}
\begin{align}\label{ABC:A}
  A_{\phi(x)} 
  &:= 1 + \nu - 3 \nu \cos^2 \phi(x),
\\\label{ABC:B}
  B_{\phi(x)} 
  &:= 2 \big( \log 2 - (1 - \log 2) \nu - (3 \log 2 - 2) \nu \cos^2 \phi(x) \big), \\\label{ABC:C}
  C_{\phi(x)} 
  &:= \frac12 \big( -3 - \nu + 3 (1 + \nu) \cos^2 \phi(x) \big)
\end{align}
\end{subequations}
and $O(\e)$ are uniformly bounded in $x$.
\end{thm}

\paragraph{Idea of the proof}  The main idea for expanding $F_\e (x)$ is to add and subtract integration over $\cC_\e(x)$. This idea dates back to \cite{GavazzaBarnett76}.
 The integral of $G \cdot \tau$ over $\cC_\e(x)$ can be expanded explicitly in terms of $\e$, which gives rise to the first two terms (i.e.\ the local contribution) in the right-hand side of \eqref{Fe:t}. The remaining two terms are the integrals over $\Gamma_\e(x)$ and $\cC_\e(x)$, which are given in the right-hand side of \eqref{Psi:formal}. The leading order terms in the expansion of these two terms cancel out, and the next order terms turn out to be (relying on Stokes' Theorem) a Cauchy sequence in $\e$, which implies the existence of the limit in \eqref{Psi:formal}. In Remark \ref{r:Psi} we specify the extension of $\tau$ from $\Gamma$ to $\cC(x)$ and the convergence rate of the limit $\e \to 0$ in \eqref{Psi:formal}.

For the expansion of $\cF_\e (x)$, we add and subtract integration over $\cC(x)$. Since $G_\e$ is regular, there is no need to remove $B(x,\e)$ from the curves. The integral over $\cC(x)$ can be expanded explicitly, which gives rise not only to the first two terms in the right-hand side of \eqref{Fe:t}, but also to the new term $\kappa(x) C_{\phi(x)}$ appearing in \eqref{cFe:t}. To expand the joint integral over $\Gamma$ and $\cC(x)$, we show that it is close in value to the integrals of the $\e$-independent integrand $G(y-x) \cdot \tau (y)$ over $\Gamma_\e(x)$ and $\cC_\e(x)$. Then, the same value $\Psi(x)$ appears naturally in the limit $\e \to 0$, and we obtain the following alternative characterization of $\Psi$:
\begin{equation} \label{Psi:formal:2}
  \Psi(x) = \lim_{\e \to 0} \oint_{\Gamma \cup \cC(x)} G_\e(y-x) \cdot \tau (y) \, dy.
\end{equation}

\paragraph{Remarks} We remark on two aspects of Theorem \ref{t:ana}.
First, instead of using $\cC_\e(x)$ in \eqref{Psi:formal}, there are many other choices to cancel out the singularity of the integrand at $y = x$. The benefit of using $\cC(x)$ is that it depends on $\Gamma$ only through the local information $x$, $\tau(x)$ and $\kappa(x)$, and that the integral of $G(y-x) \cdot \tau (y)$ over $\cC_\e(x)$ can be explicitly expanded in terms of $\e$.

Second, part of the error $O(\e)$ comes from a Taylor approximation of $\Gamma$ at $x$, which we can expand up to third order thanks to the assumption that $\Gamma$ is of class $C^3$. If $\Gamma$ would only be of class $C^{2,\alpha}$ for some $\alpha \in (0,1)$ specifying the H\"older exponent of the second derivative, then the error term becomes $O(\e^\alpha)$. 

\paragraph{Contribution} The two main contributions of Theorem \ref{t:ana} to the literature are as follows.
First, the relative error of the expansions is $O(\e)$, while the relative error in \eqref{F:expa:O1} is $O(1 / |\log \e|)$. Especially in view of the typical value $\e \sim 10^{-4}$ this is a huge improvement.

Second, the difference between the two formulas for the forces is
\begin{equation} \label{cFe:Fe:diff}  
  \cF_\e(x) - F_\e(x) = \kappa(x) C_{\phi(x)} + O(\e)
\end{equation}
as $\e \to 0$. The local nature of this difference stems from the result that the nonlocal term $\Psi(x)$ can be expressed in terms of the modeling assumptions that underly either $\cF_\e(x)$ or $F_\e(x)$; see \eqref{Psi:formal} and \eqref{Psi:formal:2}. 

The value in \eqref{cFe:Fe:diff} quantifies the difference between  two inherently different models for the self-interaction force. Since this difference is a local term which is relatively small with respect to the leading order term in \eqref{Fe:t}, there is no reason to believe that one of the two models is superior over the other. This gives a new perspective to the comparisons made in \cite[Sect.\ 3]{CaiArsenlisWeinbergerBulatov06}. In fact, in the next section we will see that the self-interaction force in both models can be discretized in a similar manner.

\subsection{Discretizations of $F_\e$ and $\cF_\e$}
\label{s:i:num}

The expansions in Theorem \ref{t:ana} of $F_\e$ and $\cF_\e$ open up new possibilities to compute them numerically from a discretized version of $\Gamma$. Such computations are necessary in discrete dislocation dynamics, which is a substantial component in the current research on plasticity. We refer to \cite{LeSarCapolungo20} for a review, and cite several specific papers below. 

The need for expansions in the numerical computation of $F_\e(x)$ and $\cF_\e(x)$ as opposed to using their definitions in \eqref{Fe} and \eqref{cFe} can be understood as follows. The curve segments of $\Gamma$ close to $x$ yield the main contribution to the integrals in \eqref{Fe} and \eqref{cFe}. Hence, to avoid large errors, one requires a fine discretization of $\Gamma$ around $x$, typically of size $O(\e)$. Since dislocation dynamics require $F_\e(x)$ and $\cF_\e(x)$ to be evaluated at many points $x$, this requires a discretization of size $O(\e)$ everywhere along $\Gamma$, which is computationally expensive. A similar reason deems phase-field models impractical too.

Such fine discretizations of $\Gamma$ are not needed when expansions such as those in Theorem \ref{t:ana} are used. Indeed, the terms in the expansions which need to be computed are independent of $\e$. Therefore, we introduce a discretization parameter $h > 0$ for the discretization of $\Gamma$, and treat it \emph{independently} of $\e$. Our second of the two aims in this paper is to develop accurate numerical schemes for $F_\e(x)$ and $\cF_\e(x)$. Our main contribution (Theorem \ref{t:num}) is to quantify the accuracy by estimating the error made when replacing $F_\e(x)$ and $\cF_\e(x)$ by their corresponding schemes.

To the best of our knowledge, there are no discretization errors beyond $O(1)$ in $\e$ available in the literature. This could explain the rather crude numerical schemes that are used in the literature. For instance, the approaches in \cite{ZbibRheeHirth98,
Schwarz99,
ArsenlisCaiTangRheeOppelstrupHommesPierceBulatov07} 
first discretize $\Gamma$ to a polygon, and then compute the self-interaction force from the known formula for straight dislocation segments (see \eqref{Ixy} below). A special treatment for connecting segments is made to avoid the singularity in $G$. It is difficult to track the discretization error made in this manner. In fact, our approach below reveals that an approximation by straight lines may cause a large error if not dealt with carefully. Another set of examples are the schemes used in \cite{GhoniemTongSun00,
ZhuChapmanAcharya2013} and in the paper series by Bene\v s, Kratochv\'il, Pau\v s \textit{et al.} (see \cite{KolarBenesKratochvilPaus18} and references therein), which are based on the expansion in \eqref{F:expa:O1}, but neglect the nonlocal contribution $\Psi(x)$. This creates a relative discretization error of size $O(1/|\log \e|)$. 

\begin{figure}[ht]
\centering
\begin{tikzpicture}[scale=1.5, >= latex]
\def \rr {.04} 

\draw (0,0) ellipse (2 and 1);

\fill (2,0) circle (\rr);
\fill (-2,0) circle (\rr);
\fill (0,1) circle (\rr);
\fill (0,-1) circle (\rr);

\foreach \x in {0,180}{ 
\foreach \pmm in {1,-1}{ 
  \begin{scope}[rotate = \x]  
    \fill (\pmm*1.732,.5) circle (\rr);
    \fill (\pmm*1,.866) circle (\rr);  
    \draw[red] (\pmm*2,0) -- (\pmm*1.732,.5) -- (\pmm*1,.866) -- (0,1);   
  \end{scope}
}
}

\draw (2,0) node[right] {$x_1 = x_{N+1}$};
\draw (1.732,.5) node[anchor = south west] {$x_2$};
\draw (1,.866) node[above] {$x_3$};
\draw (1.732,-.5) node[anchor = north west] {$x_{N} = x_0$};
\draw (1,-.866) node[below] {$x_{N-1}$};
\draw (-1.732,.5) node[anchor = south east] {$x_i$};
\draw (-1,.866) node[above] {$x_{i-1}$};
\draw[red] (-1.3,.68) node[below] {$\gamma_i$};
\draw (-1.4,-.78) node[anchor = north east] {$\Gamma$};
\draw[red] (-.5,-.93) node[above] {$\Gamma^h$};
\end{tikzpicture} \\
\caption{Sketch of a possible discretization $\Gamma^h$ of $\Gamma$.}
\label{fig:Gammah}
\end{figure}
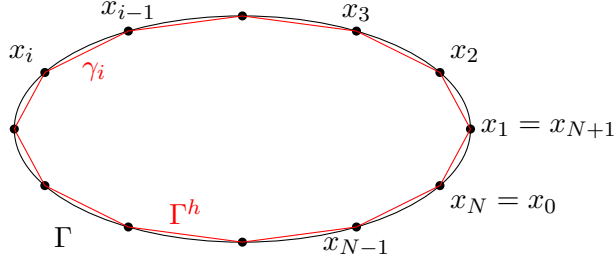

Our second main result, Theorem \ref{t:num}, guarantees a much smaller discretization error. To define the corresponding numerical schemes, we first discretize $\Gamma$. Let $h > 0$ be a spatial discretization parameter, which we assume to be small enough with respect to $\Gamma$. Let $x_i \in \Gamma$ be discretization points for $i = 1,\dots,N$. Figure \ref{fig:Gammah} illustrates this setting. We consider $\bx := (x_1, \ldots, x_N) \in (\R^2)^N$ as the complete list of variables which describe the discretization of $\Gamma$. Here and henceforth, we extend $x_i$ periodically over its index by setting $x_{i + Nj} := x_i$ for any $j \in \Z$. In addition, we choose the ordering of $x_1, \ldots, x_N$ such that, when transversing $\Gamma$ in counter-clockwise direction, the points $x_i$ appear with increasing index.
We set
\begin{equation*}
  \gamma_i := \{ t x_i + (1-t) x_{i-1} : 0 \leq t \leq 1 \}
  \qquad \text{for } i = 1,\dots,N
\end{equation*}
as the closed line segments connecting $x_{i-1}$ to $x_i$, and set 
\begin{equation*} 
  \Gamma^h := \bigcup_{i=1}^N \gamma_i
\end{equation*}
as the piecewise-affine closed curve which discretizes $\Gamma$.
To relate the choice of $\bx$ to $h$, we assume that $|\gamma_i| \sim h$, i.e.\ there exists a universal constant $C > 1$ such that
\begin{equation} \label{h:xi}
  \frac1C \max_{1 \leq i \leq N}|\gamma_i| \leq h \leq C \min_{1 \leq i \leq N}|\gamma_i|.
\end{equation}

Our discretization of $F_\e(x_i)$ is
\begin{equation} \label{Feih}
  F_{\e,i}^h (\bx)
  :=  \frac{ \kappa^h(x_i) A_{\phi^h(x_i)}}2 \log \frac{|x_{i+m^h} - x_i| |x_{i-m^h} - x_i|}{\e^2} + \sum_{j = m^h + 1}^{N - m^h} \int_{\gamma_{i + j}} G(y-x_i) \cdot \tau (y) \, dy,
\end{equation}
where
\begin{equation} \label{mh}
  m^h := \big\lceil h^{-1/3} \big\rceil \in \N,
\end{equation}
$\kappa^h(x_i)$ is the approximation of $\kappa(x_i)$ given by
\begin{equation} \label{kaphxi}
  \kappa^h(x_i) 
  := 2 \Big( \frac1{|y_+|} + \frac1{|y_-|} \Big) \frac{ y_- \cdot Q y_+ }{ (|y_+| + |y_-|)^2 }
\end{equation}
with
\begin{equation} \label{ypm}  
  y_\pm := x_{i \pm 1} - x_i
  \qand
  Q := \begin{bmatrix} 0 & 1 \\ -1 & 0 \end{bmatrix},
\end{equation}
and $A_{\phi^h(x_i)}$ is explicit as a function of $\cos \phi^h(x_i) = n^h(x_i) \cdot e_1$ (see \eqref{ABC:A}), where 
\begin{equation} \label{nhxi} 
  n^h(x_i) 
  := \frac{\tilde n^h(x_i)}{|\tilde n^h(x_i)|} \qtext{with}
  \tilde n^h(x_i) 
  := Q \Big( |y_-| \frac{y_+}{|y_+|} - |y_+| \frac{y_-}{|y_-|} \Big)
\end{equation}
is an approximation of $n(x_i)$. The expression \eqref{Feih} resembles the definition of $F_\e$ in \eqref{Fe} rather than its expansion in \eqref{Fe:t}. Indeed, if we replace in \eqref{Fe} the curve $\Gamma_\e$ by the discretized curve $\Gamma^h$ and remove from it several line segments $\gamma_j$ close to the point $x_i$, we obtain \eqref{Feih} with a local correction term. Yet, the construction of \eqref{Feih} relies on Theorem \ref{t:ana}. 

Finally, we introduce our discretization of $\cF_\e(x_i)$. It is simply given by
\begin{equation} \label{cFeih}
  \cF_{\e,i}^h (\bx)
  := F_{\e,i}^h (\bx)
    + \kappa^h(x_i) C_{\phi^h(x_i)},
\end{equation}
where $C_{\phi^h(x_i)}$ is -- similar to $A_{\phi^h(x_i)}$ -- a function of the two components of the vector $n^h(x_i)$.

\begin{thm}[Discretization of $F_\e(x)$ and $\cF_\e(x)$] \label{t:num}
Let $\Gamma$ be a non-intersecting closed curve of class $C^3$ with finite length. Let $\e, h \in (0,1)$ be small enough with respect to $\Gamma$, and let $\bx$ be as above such that \eqref{h:xi} holds.
Then, for all $i = 1,\ldots,N$ 
\begin{align} \label{err:t}
  \big| F_{\e,i}^h (\bx) - F_\e(x_i) \big| + \big| \cF_{\e,i}^h (\bx) - \cF_\e(x_i) \big| 
  \leq C (\e + h |\log \e| + h^{2/3}),
\end{align}
where $C > 0$ is independent of $\e, h, \bx$.
\end{thm}

\paragraph{Idea of the proof} 
The proof is constructive, and motivates the scheme in \eqref{Feih} by starting from the expansion of $F_\e(x_i)$ in \eqref{Fe:t}. This explains the appearance of $\e$ in the error term. Since we construct $\kappa^h(x_i)$ to be $O(h)$ close to $\kappa(x_i)$, the logarithmic term in \eqref{Fe:t} leads to the error term $h |\log \e|$ in \eqref{err:t}.

The interesting part is the discretization of $\Psi(x_i)$. The characterization of $\Psi(x_i)$ in \eqref{Psi:formal} is constructed carefully for the singularity in the integrand $G(y-x_i) \cdot \tau(y)$ to cancel out. By adding and subtracting integration over $\cC(x_i)$ in this characterization, we construct effectively an approximation of $\Gamma$ around $x_i$ up to third order. However, the approximation of $\Gamma$ by the polygon $\Gamma^h$ is only of first order, which is not enough to cancel out the singularity. Hence, the replacement of $\Gamma_\e(x_i) = \Gamma \setminus B(x_i,\e)$ in the definition of $\Psi(x_i)$ in \eqref{Psi:formal} by $\Gamma^h \setminus B(x_i,\e)$ is expected to yield an error which is large as $\e \to 0$.

To work around this problem, we rely on a by-product from the proof of Theorem \ref{t:ana} (see Remark \ref{r:Psi}). This by-product states that the integral in \eqref{Psi:formal} is $O(\e)$ close to $\Psi(x_i)$. This creates a trade-off with the previously mentioned error from replacing  $\Gamma_\e (x_i)$ by $\Gamma^h \setminus B(x_i,\e)$, which becomes larger as $\e$ tends to $0$. Instead of applying these error estimates for the given value of $\e$ in Theorem \ref{t:num}, we apply them for a possibly different, $h$-dependent 
$\e^h$ for which both error terms are of the same order. It will turn out that
\[
  \e^h \sim h^{2/3}.
\]
This explains the appearance of $h^{2/3}$ in the error in Theorem \ref{t:num}. It also explains the use of the number $m^h$ in \eqref{mh}, which is chosen such that
\[
  |x_{i \pm m^h} - x_i| \sim \e^h.
\]

\paragraph{Remarks}
We remark on five aspects of Theorem \ref{t:num}.
First, note that Theorem \ref{t:num} requires no estimates between $\e$ and $h$, which means that \eqref{err:t} is valid in the two-dimensional parameter space $\e, h > 0$ in a (complete) neighborhood around $0$. For practical purposes, $\e > 0$ is a modeling parameter, and the choice of $h$ is up to the discretion of the user. Then, the following corollary of Theorem \ref{t:num} is of more practical use: given the assumptions of Theorem \ref{t:num}, it holds that
\begin{equation} 
  \e^{3/2} \leq h \leq |\log \e|^{-3}
  \implies
  \big| F_{\e,i}^h (\bx) - F_\e(x_i) \big| + \big| \cF_{\e,i}^h (\bx) - \cF_\e(x_i) \big| 
  \leq C h^{2/3}
\end{equation}
for some constant $C > 0$ independent of $\e, h, \bx$.

Second, the statement ``$\e,h$ are small enough with respect to $\Gamma$" can be made more precise. In fact, the proof requires $\e$ and $h$ to be small enough with respect to the four constants
\begin{equation} \label{Gamma:glo:csts}
  |\Gamma|, \quad
  \max_{s \in \R} |\varphi''(s)|, \quad
  \max_{s \in \R} |\varphi'''(s)|, \quad
  \min_{s<t} \frac{t-s}{|\varphi(t) - \varphi(s)|},
\end{equation}
where $|\Gamma|$ is the length of $\Gamma$ and $\varphi$ is an arc length parametrization of $\Gamma$. The fourth constant above is large if, for instance, in the situation on the right in Figure \ref{fig:Gamma0}, a second part of $\Gamma$ would be inside $B(x,\e)$. 

Third, note the practical property that $F_{\e,i}^h$ and $\cF_{\e,i}^h$ can be directly computed from $\bx$. Indeed, the first term in \eqref{Feih} and the second term in \eqref{cFeih} are explicitly expressed in terms of the five points $x_{i-m^h}$, $x_{i-1}$, $x_i$, $x_{i+1}$ and $x_{i+m^h}$. To express the remaining second term in \eqref{Feih} explicitly in terms of $\bx$, we recall (see, e.g.\ \cite{DeWit60,
HirthLothe82}) 
the following formula for the force that a line segment $\gamma_{x \to y}$ which connects $x \in \R^2$ to $y \in \R^2$ exerts on a curve segment at $0$ in normal direction:
\begin{equation} \label{Ixy} 
  I(x,y) :=
  \int_{\gamma_{x \to y}} G(z) \cdot \tau(z) \, dz
  = \frac1{ |x| |y| + x \cdot y }
    \Big(  \frac x{|x|} + \frac y{|y|} \Big)^T
    \begin{bmatrix}
  0 & -1 \\
  1 - \nu & 0
  \end{bmatrix} (y-x). 
\end{equation}
Using this explicit formula, the second term in \eqref{Feih} becomes
\[
  \sum_{j = m^h + 1}^{N - m^h} \int_{\gamma_{i + j}} G(y-x_i) \cdot \tau (y) \, dy
  = \sum_{j = m^h + 1}^{N - m^h} I(x_{i+j-1} - x_i, x_{i+j} - x_i),
\]
which is indeed an explicit function of $x_i$ and $x_{m^h + i}, x_{m^h + i + 1}, \ldots, x_{N-m^h + i}$.

Fourth, the discretizations $F_{\e,i}^h(\bx)$ and $\cF_{\e,i}^h(\bx)$ are numerically suboptimal, because they do not depend on $x_{i \pm j}$ for $j = 2,3,\ldots,m^h-1$. One could decide to require higher regularity on $\Gamma$, and use these points $x_{i \pm j}$ to develop a higher order approximation for $n(x)$ and $\kappa(x)$ than those in \eqref{nhxi} and \eqref{kaphxi}. This alone, however, will not decrease the size of the discretization error in \eqref{err:t}.

Fifth, the description of the proof reveals that the relatively large error term $h^{2/3}$ appears because of the low regularity (Lipschitz, but not $C^1$) of the discretized curve $\Gamma^h$. By using a higher order discretization based on splines (see e.g.\ \cite{GhoniemTongSun00}) it might be possible to get this part of the error down to $O(h)$, simply by taking $m^h = 1$ and $\e^h \sim h$. We leave this direction as an open problem.

\paragraph{Contribution} The main contribution of Theorem \ref{t:num} to the literature is that, for the first time, it introduces and estimates the discretization error made when computing the self-interaction force from a numerical scheme. We consider this an important step towards building more accurate methods for computing discrete dislocation dynamics. 

\subsection{Conclusion}

Rather than the statements of Theorems \ref{t:ana} and \ref{t:num} we consider their proofs to be the main contribution of this paper. Indeed, these two theorems are merely stated for the simplest setting of a two-dimensional, isotropic medium, and they only entail two out of several models for the self-interaction force. Yet, we expect the methodology (i.e.\ the skeleton of the proof) to apply to different models for this force (especially those constructed by smearing out the dislocation by a convolution kernel) and to a three-dimensional setting. 

Finally, we recall that expansions such as those in Theorem \ref{t:ana} account for the nonlocal linear-elastic contribution of the self-interaction force of a dislocation, but that they do not provide a proper accounting for the contribution from the dislocation core. Yet, expansions as those in Theorem \ref{t:ana} reveal the connection between different models for the nonlocal contribution of the self-interaction force, which gives more freedom in constructing an appropriate coupling with atomistic models. 

The remainder of the paper is concerned with proving Theorems \ref{t:ana} and \ref{t:num}. Section \ref{s:t:ana} contains the proof of Theorem \ref{t:ana} and Section \ref{s:t:num} contains the proof of Theorem \ref{t:num}.

\section{Proof of Theorem \ref{t:ana}}
\label{s:t:ana} 

The notation convention in Sections \ref{s:t:ana} and \ref{s:t:num} is as follows. Whenever convenient, we denote balls such as $B(x,\e)$ by $B_\e(x)$. We reserve $c, C$ for generic positive constants which do not depend on any of the important variables. We use $C$ in upper bounds (and think of it as possibly large) and $c$ in lower bounds (and think of it as possibly small). While $c, C$ may vary from line to line, in the same display they refer to the same value. If different constants appear in the same display, we denote them by $C, C', C'', \ldots$. Finally, in local computations we sometimes reuse notation. For instance, in the computation of line integrals, the parametrization of the corresponding curve is an auxiliary step which is of no further use outside of the computation of the line integral; we use the symbol $\varphi$ to denote various parametrizations. 

We prove Theorem \ref{t:ana} by establishing \eqref{Fe:t} and \eqref{cFe:t} in respectively Sections \ref{s:t:ana:1} and \ref{s:t:ana:2}.

\subsection{Expansion of $F_\e$}
\label{s:t:ana:1}

In this section we construct the expansion of $F_\e(x)$ as given in \eqref{Fe:t}, which is the first of the two parts of Theorem \ref{t:ana}. We fix an $x \in \Gamma$, translate the coordinate system such that $x$ is at the origin, and remove $x$ from the notation. Then, the definition of $F_\e(x)$ in \eqref{Fe} reads as
\[
  F_\e
  = \int_{\Gamma_\varepsilon} G(y) \cdot \tau (y) \, dy
  = \int_{\Gamma_\varepsilon} G \cdot \tau.
\]
Regarding \eqref{ktaun0}, we express the curvature, tangent vector (counter-clockwise direction) and outward pointing normal vector of $\Gamma$ at $0$ respectively by $\kappa_0$, $\tau_0$ and $n_0$ (to avoid clutter, we denote the related angle $\phi_0$ simply by $\phi$). The statement in Theorem \ref{t:ana} that the error term is uniform in $x$ translates to the requirement that the error term $O(\e)$ has to be uniform in $\phi \in [0, 2\pi)$, $\kappa_0$ and other local information of $\Gamma$ at $0$.
Depending on the sign of $\kappa_0$, we split three cases.

\paragraph{The case $\kappa_0 > 0$.} Let $r_0 = 1/\kappa_0$ be the radius of the tangent circle $\cC = \partial B( r_0 n_0, r_0)$ of $\Gamma$ at $0$. For technical reasons, we first assume that $\Gamma \cap \cC = \{0\}$, i.e., $\Gamma$ and $\cC$ intersect only at $0$. We treat the general case afterwards. 

The idea for expanding $F_\e$ is to replace $\Gamma_\e$ by
$$
  \cC_\varepsilon := \cC \setminus B_\varepsilon(0),
$$ 
and to show that the error made by this replacement is of the form $\Psi + O(\e)$. With this aim, let $\e$ be small enough such that $\partial B_\varepsilon(0)$ intersects with $\Gamma$ and $\cC$ at precisely two points. Note that this upper bound on $\e$ can be constructed from the constants in \eqref{Gamma:glo:csts}, such that it does not depend on the local information of $\Gamma$ at $0$. Let $\gamma_\e \subset \partial B_\varepsilon(0)$ be the union of the two disjoint arcs which connect the endpoints of $\Gamma_\e$ and $\cC_\e$. Figure \ref{fig:Gammae:Ce} illustrates these curves. We extend $\tau$ (i.e., the direction of the tangent vector on $\Gamma$) to $\cC_\e$ and $\gamma_\e$ such that there is a consistent direction in which the closed curve $\Gamma_\e \cup \cC_\e \cup \gamma_\e$ is traversed. In particular, this means that $\tau$ is such that $\cC_\e$ is traversed in counter-clockwise direction. Then, we decompose 
\begin{equation} \label{Ie123}
  F_\e 
  = \oint_{\Gamma_\e \cup \cC_\e \cup \gamma_\e} G \cdot \tau
    - \int_{\cC_\varepsilon} G \cdot \tau
    - \int_{\gamma_\varepsilon} G \cdot \tau
  =: F_\e^1 - F_\e^2 - F_\e^3,
\end{equation}
and expand $F_\e^1$, $F_\e^2$ and $F_\e^3$ independently. 

\begin{figure}[ht]
\centering
\begin{tikzpicture}[scale=1.2, >= latex]
\def \rr {.04} 

\begin{scope}[scale = 1.4] 
\draw [thick] plot [smooth cycle, tension=.7] coordinates {(0,0) (-.1,1) (.8,1) (1.4,-.5) (1, -1.2) (-.1, -1)};
\end{scope}
\draw [thick, blue] (-1.5,0) circle (1.5);

\fill [white] (0,0) circle (1);
\draw[green!70!black, dotted] (0,0) circle (1);


\begin{scope}[scale = 1.4] 
\draw [dotted] plot [smooth cycle, tension=.7] coordinates {(0,0) (-.1,1) (.8,1) (1.4,-.5) (1, -1.2) (-.1, -1)};
\draw (1.4,-.5) node[right]{$\Gamma_\e$};
\end{scope}
\draw [dotted, blue] (-1.5,0) circle (1.5);

\begin{scope}[rotate = 101] 
  \draw[green!70!black, thick] (1,0) arc (0:9:1);       
\end{scope}
\begin{scope}[rotate = 250] 
  \draw[green!70!black, thick] (1,0) arc (0:8:1);       
\end{scope}

\draw[dotted] (0,0) -- (-1.5,0);
\begin{scope}[shift={(-1.5,0)}, rotate = 30] 
  \draw[blue, thin] (0,0) -- (0,1.5) node[midway, left]{$r_0$};       
\end{scope}
\draw[blue] (-3,0) node[left]{$\cC_\e$};
\draw[red, <->] (0,-.8) node[right]{$\tau_0$} -- (0,0) -- (-.8,0) node[below]{$n_0$};
\draw[green!70!black] (-.8,-.6) node[left]{$B_\e(0)$};
\draw[green!70!black] (1,-1) node {$\gamma_\e$};
\begin{scope}[shift={(1,-1)}, rotate = 180] 
  \draw[->, green!70!black] (0.2,0) -- (1.1,0);    
\end{scope}
\begin{scope}[shift={(1,-1)}, rotate = 123] 
  \draw[->, green!70!black] (0.2,0) -- (2.2,0);    
\end{scope}
\draw (1,1) node {$\Omega$};

\fill (0,0) circle (\rr);
\fill[blue] (-1.5,0) circle (\rr);
\end{tikzpicture} \\
\caption{Sketch of the closed curve formed by $\Gamma_\e$, $\cC_\e$ and $\gamma_\e$.}
\label{fig:Gammae:Ce}
\end{figure} 

We start with $F_\e^3$. Let 
$$
  \varphi (\theta) = \varepsilon \begin{bmatrix} \cos \theta \\ \sin \theta \end{bmatrix} \qquad \text{with } \alpha < \theta < \alpha + \delta_\varepsilon 
$$ 
be a parametrization of one of the two arcs of $\gamma_\varepsilon$, where $\alpha = \phi - \tfrac\pi2 + O(\varepsilon)$ (recall $\phi$ from \eqref{ktaun0}). For $\e$ small enough, it follows from the $C^3$-regularity of $\Gamma$ that the endpoints of $\cC_\e$ and $\Gamma_\e$ are a distance of $O(\varepsilon^3)$ apart. Hence, $\delta_\varepsilon = O(\varepsilon^2)$. Then, the contribution of this arc to the value of $F_\e^3$ is
\begin{equation*} 
  \int_{\alpha}^{\alpha + \delta_\varepsilon} 
      G(\varphi(\theta)) \varphi'(\theta) \, d\theta \\
  = -\frac1\varepsilon \int_{\alpha}^{\alpha + \delta_\varepsilon} 
      \big( (1-\nu) \sin^2 \theta + \cos^2 \theta \big) \, d\theta
  = O(\delta_\varepsilon/\varepsilon)  
  = O(\varepsilon).
\end{equation*}
The contribution of the second arc of $\gamma_\e$ can be treated analogously. This proves that $F_\e^3 = O(\varepsilon)$.

Next we expand the term $F_\e^1$ in \eqref{Ie123}.  With this aim, we first show that $(F_\e^1)_\e$ is a Cauchy sequence. Let $\e > 0$ be small enough, and take $\delta \in (0, \e)$. Let $\Omega$ be the region enclosed by $\Gamma$. Then, 
\[
  F_\e^1 - F_\delta^1
  = \oint_{\partial \omega_{\delta, \e}} G \cdot \tau,
\]
where the open set
\begin{equation} \label{omde}
  \omega_{\delta, \e} 
  := B_\e(0) \setminus \overline{ B_\delta(0) \cup B_{r_0}( r_0 n_0) \cup \Omega } 
\end{equation}
is a subset of the narrow wedges between $\Gamma$ and $\cC$ (see Figure \ref{fig:Gammae:Ce}). Since $\Gamma$ and $\cC$ intersect only at $0$, $\partial \omega_{\delta, \e}$ consists of two disjoint closed loops. Then, by Stokes' Theorem,
\[
  \oint_{\partial \omega_{\delta, \e}} G \cdot \tau  
  = \iint_{\omega_{\delta, \e}} g,
\]
where
\begin{equation} \label{g} 
  g(y) := \curl G(y)
  = \frac{ (1+\nu) y_1^2 + (1 - 2\nu) y_2^2 }{|y|^5}.
\end{equation}
Hence,
\begin{equation} \label{omde:int:est}
  \bigg| \oint_{\partial \omega_{\delta, \e}} G \cdot \tau \bigg|  
  \leq C \iint_{\omega_{\delta, \e}} \frac1{ |y|^3 } \, dy
\end{equation}
for some constant $C$ which only depends on $\nu$. To estimate this integral, we use polar coordinates to write
\begin{equation} \label{OmOme:pmzd}
  \omega_{\delta, \e}
  = \{ (r, \theta) : \delta < r < \varepsilon, \, \theta \in \Theta(r) \}, 
\end{equation} 
where $\Theta(r) \subset \R / (2\pi \Z)$ is the union of two intervals. Analogously to the argument for $|\gamma_\e| = O(\e^3)$, we obtain that $|\Theta(r)| = O(r^2)$. Then, 
\begin{equation} \label{OmOme:intd}
  \iint_{\omega_{\delta, \e}} \frac1{ |y|^3 } \, dy
  = \int_\delta^\varepsilon \int_{\Theta(r)} \frac1{r^3} r \, d\theta dr
  = \int_\delta^\varepsilon \frac1{r^2} |\Theta(r)| \, dr
  = O(\varepsilon)
\end{equation}
uniformly in $\delta$.
Hence, $(F_\e^1)_\e$ is a Cauchy sequence, and 
\begin{equation*} 
  F_\e^1 
  = \Psi + O(\varepsilon)
\end{equation*}
for some $\Psi \in \R$.

Finally, we expand the term $F_\e^2$ in \eqref{Ie123}. With this aim, we parametrize the arc $\cC_\varepsilon$ in  counter-clockwise direction by
\begin{align} \label{Ie2:vphi}
  \varphi (\theta) := r_0 \begin{bmatrix} \cos \phi + \cos (\theta + \phi + \pi) \\ \sin \phi + \sin (\theta + \phi + \pi) \end{bmatrix} \qquad \text{with } \alpha < \theta < 2\pi - \alpha,
\end{align}
where $\alpha \in (0, \pi)$ is such that 
$$
  \varepsilon^2 = |\varphi(\alpha)|^2.
$$
For simplicity, we set
\begin{align*}
  c_\phi &:= \cos \phi & c &:= \cos ( \tfrac\theta2 ) \\
  s_\phi &:= \sin \phi & s &:= \sin ( \tfrac\theta2 )
\end{align*}
and 
\[
  Q_\phi := \begin{bmatrix} c_\phi & -s_\phi \\ s_\phi & c_\phi \end{bmatrix}
\]
as the counter-clockwise rotation matrix by angle $\phi$. Using trigonometric identities, we get
\begin{align} \label{Ie2:vphi:props}
  \varphi (\theta) = 2 r_0 s Q_\phi \begin{bmatrix} s \\ -c \end{bmatrix}, \quad
  | \varphi (\theta) | = 2 r_0 |s| = 2 r_0 s
  \quad \text{and} \quad
  \varphi' (\theta) = r_0 Q_\phi \begin{bmatrix} 2sc \\ s^2 - c^2   \end{bmatrix}.
\end{align}
In particular,
\[
  \e^2 = |\varphi(\alpha)|^2 = 4 r_0^2 \sin^2 (\alpha/2).
\]
Hence,
\begin{equation*} 
  \alpha = \frac{\e}{r_0} + O(\e^3),
\end{equation*}
where $O(\e^3)$ can be bounded uniformly in $\kappa_0$ since $1/r_0 = \kappa_0 \leq \max_\Gamma |\kappa|$.

Noting that $s > 0$, we obtain from the preparations above that
\begin{align} \notag 
  F_\e^2 
  &= \int_{\alpha}^{2\pi - \alpha} 
      G(\varphi(\theta))\varphi'(\theta) \, d\theta \\\notag
  &= \int_{\alpha}^{2\pi - \alpha} 
      \frac{1}{4 r_0 s^2} 
      \begin{bmatrix} s & -c \end{bmatrix}
      Q_\phi^T
      \begin{bmatrix}
  0 & -1 \\
  1-\nu & 0
  \end{bmatrix} 
  Q_\phi
      \begin{bmatrix} 2sc \\ s^2 - c^2 \end{bmatrix} 
      \, d\theta \\\label{Ie2:comp2}
  &= \frac{\kappa_0}4 \sum_{k=0}^3 C_k \int_{\alpha}^{2\pi - \alpha} 
      c^k s^{1-k} \, d\theta,
\end{align}
for some explicit constants $C_k$ which only depend on $\nu$ and $\phi$. By the symmetries of the sine and cosine, the integrands corresponding to odd powers of the cosine (i.e., $k = 1, 3$) yield zero integral. Then, substituting the explicit values of $C_0$ and $C_2$ and using $c^2 = 1 - s^2$, we obtain
\begin{equation} \label{I2e:1}
  F_\e^2 
  = \frac{\kappa_0}4 \bigg( 2 \nu (1 - 2 c_\phi^2) \int_{\alpha}^{2\pi - \alpha} 
      \sin ( \tfrac\theta2 ) \, d\theta + (3 \nu c_\phi^2 - \nu - 1) \int_{\alpha}^{2\pi - \alpha} 
      \frac1{\sin ( \tfrac\theta2 )} \, d\theta \bigg).
\end{equation}

Computing the integrals, we get
\begin{align*}
  \int_{\alpha}^{2\pi - \alpha} \sin ( \tfrac\theta2 ) \, d\theta
  = 4 + O(\alpha^2)
  = 4 + O(\e^2)
\end{align*}
and
\begin{multline*}
  \int_{\alpha}^{2\pi - \alpha} 
      \frac1{\sin ( \tfrac\theta2 )} \, d\theta
  = 2 \log \Big( \tan \frac\theta4 \Big) \bigg|_{\theta = \alpha}^{2\pi - \alpha} 
  = 4 \log \frac1{\alpha} + 8 \log 2 + O(\alpha^2)  \\
  = 4 \log \frac{r_0}\e + 8 \log 2 + O(\e^2). 
\end{multline*}
Collecting the computations above, we obtain
\begin{equation} \label{Ie2:expanded}
  F_\e^2 
  = -\kappa_0 \Big( A_\phi \log \frac1{\varepsilon \kappa_0} + B_\phi + O (\e^2) \Big),    
\end{equation} 
where $A_\phi$ and $B_\phi$ are defined in \eqref{ABC}.

Finally, collecting all estimates on $F_\e^i$, \eqref{Fe:t} follows for the case in which $\Gamma \cap \cC = \{0\}$.
\medskip

Next we generalize the proof of \eqref{Fe:t} to general curves $\Gamma$ without any restrictions to the number of intersections with $\cC$. In this case, the derivation of the expansions of $F_\e^2$ and $F_\e^3$ above remain valid. The derivation of the expansions of $F_\e^1$ requires minor modifications. Indeed, since $\Omega$ and $B_{r_0} (r_0 n_0)$ may overlap, \eqref{omde} may not be the correct region to consider.  

To find an alternative definition for $\omega_{\delta, \e}$, we take $\e$ small enough so that both $\Gamma \cap B_\e(0)$ and $\cC \cap B_\e(0)$ can be parametrized by height functions $h_1, h_2$ on the part $\{t \tau_0 : |t|  < \e \}$ of the tangent line of $\Gamma$ at $0$ (one can use Figure \ref{fig:Gammae:Ce} for a visualization). Note that
\begin{align*}
  \Gamma \cap B_\e(0) &\subset \left\{ t \tau_0 + h_1(t) n_0 : -\e < t < \e \right\} \\
  \cC \cap B_\e(0) &\subset \left\{ t \tau_0 + h_2(t) n_0 : -\e < t < \e \right\}.
\end{align*}
Then, we define $\omega_{\delta,\e}$  as the open set inside the annulus $B_\e(0) \setminus B_\delta(0)$ which is in between $\Gamma$ and $\cC$. Note that this definition is consistent with that in \eqref{omde}. The connected components of the closed set
\[
  \{ t \in [-\e, \e] : h_1(t) = h_2(t) \}
\]
separate $\omega_{\delta,\e}$ into connected components. If the number of components is finite, then we can apply Stokes' Theorem as in \eqref{omde:int:est} on each of the components, and apply the estimate in \eqref{OmOme:intd} to conclude that $F_\e^1 = \Psi + O(\e)$. If the number of components of $\omega_{\delta,\e}$ is infinite, then we order them in size, observe that this ordering shows that the number of components is countable, write
\[
  \omega_{\delta,\e} = \bigcup_{k=1}^\infty \omega_k
\]
with $|\omega_1| \geq |\omega_2| \geq \ldots$, and use that $\{\omega_k\}_{k=1}^\infty$ are disjoint to estimate
\[
  \bigg| \oint_{\partial \omega_{\delta, \e}} G \cdot \tau \bigg|
  \leq \sum_{k=1}^\infty \bigg| \oint_{\partial \omega_k} G \cdot \tau \bigg|
  = \sum_{k=1}^\infty \bigg| \iint_{\omega_k} g \bigg|
  \leq \iint_{\omega_{\delta, \e}} | g |.
\]
The remainder of the argument for the expansion of $F_\e^1$ is analogous to the previous case in which $\Gamma \cap \cC = \{0\}$. This completes the proof of \eqref{Fe:t} for the case $\kappa_0 > 0$.
\medskip

\paragraph{The case $\kappa_0 = 0$.}  
This case can be treated with two minor modifications to the proof in the case $\kappa_0 > 0$. We list these modifications below.

The first modification is that we put $\cC$ as the tangent line of $\Gamma$ at $0$. Rather than a modification, this choice of $\cC$ is the natural extension of the tangent circle as $\kappa_0 \downarrow 0$. While $\cC$ is unbounded, it follows from the quadratic decay of $G(z)$ as $|z| \to \infty$ that the splitting in \eqref{Ie123} remains valid.

The second modification is the expansion of $F_\e^2$. Since $G$ is odd, we get immediately
\[
  F_\e^2 = \int_{\cC_\e} G \cdot \tau = 0,
\]
which holds without any error terms depending on $\e$.

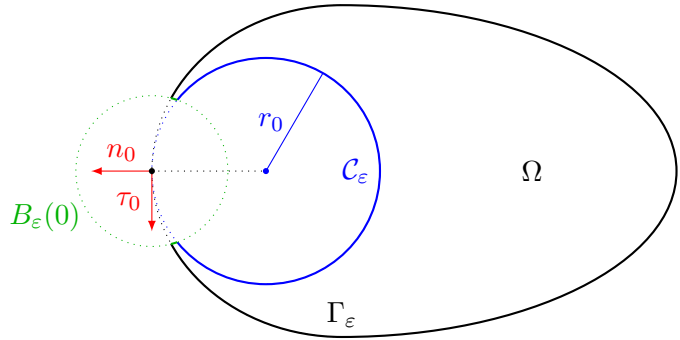
\begin{figure}[ht]
\centering
\begin{tikzpicture}[scale=1, >= latex]
\def \rr {.04} 

\draw [thick] (2.5,0) ellipse (2.5 and 2.2);
\draw [thick, blue] (1.5,0) circle (1.5);

\fill [white] (0,0) circle (1);
\draw[green!70!black, dotted] (0,0) circle (1);


\draw [dotted] (2.5,0) ellipse (2.5 and 2.2);
\fill [white] (3.1, -2.3) rectangle (5.1,2.3);
\fill [white] (2.5, 1.6) rectangle (5.1,2.3);
\fill [white] (2.5, -2.3) rectangle (5.1,-1.6);
\draw (2.5,-2.2) node[above]{$\Gamma_\e$};
\draw [dotted, blue] (1.5,0) circle (1.5);

\begin{scope}[shift={(2.5,0)}, xscale = 2] 
  \draw[thick] (0,-2.2) arc (-90:90:2.2);       
\end{scope}

\begin{scope}[rotate = 70.5] 
  \draw[green!70!black, thick] (1,0) arc (0:5:1);       
\end{scope}
\begin{scope}[rotate = 284.5] 
  \draw[green!70!black, thick] (1,0) arc (0:5:1);       
\end{scope}

\draw[dotted] (0,0) -- (1.5,0);
\begin{scope}[shift={(1.5,0)}, rotate = -30] 
  \draw[blue, thin] (0,0) -- (0,1.5) node[midway, left]{$r_0$};       
\end{scope}
\draw[blue] (3,0) node[left]{$ \cC_\e$};
\draw[red, <->] (0,-.8) -- (0,0) node[midway,left]{$\tau_0$} -- (-.8,0) node[midway,above]{$n_0$};
\draw[green!70!black] (-.8,-.6) node[left]{$B_\e(0)$};
\draw (5,0) node {$\Omega$};

\fill (0,0) circle (\rr);
\fill[blue] (1.5,0) circle (\rr);
\end{tikzpicture} \\
\caption{Sketch of a situation as in Figure \ref{fig:Gammae:Ce} when $ \kappa_0 < 0$.}
\label{fig:Gammae:kneg:Ce}
\end{figure}

\paragraph{The case $\kappa_0 < 0$.} Let $r_0 = 1/|\kappa_0|$ be the radius of the tangent circle $\cC = \partial B_{r_0}( -r_0 n_0)$.  
As in the case $\kappa_0 > 0$, we extend $\tau$ from $\Gamma$ to $ \cC_\e$ and $ \gamma_\e$ such that $\Gamma_\e \cup \cC_\e \cup \gamma_\e$ has a consistent direction in which it is traversed. In particular, $\cC_\e$ is traversed in clockwise direction, which is opposite to the case $\kappa_0 > 0$. See Figure \ref{fig:Gammae:kneg:Ce} for a sketch in the simple situation where $ \cC \subset \Omega \cup \{0\}$.  

The splitting of $F_\e$ in \eqref{Ie123} in terms of $F_\e^1, F_\e^2, F_\e^3$ still holds, and similar to the case $\kappa_0 > 0$ we obtain $ F_\e^1 = \Psi + O(\e)$ and $ F_\e^3 = O(\e)$. To compute $F_\e^2$, we parametrize $\cC_\e$ by $\tilde \varphi (\theta) := - \varphi (\theta)$, where $\varphi$ is the parametrization used in the case $\kappa_0 > 0$; see \eqref{Ie2:vphi}. Since this parametrization traverses $\cC_\e$ in the opposite direction (i.e.\ counter-clockwise), we obtain as in \eqref{Ie2:comp2}
\begin{align*}
  F_\e^2 
  &= \int_{\cC_\varepsilon} G \cdot \tau \\
  &= -\int_{\alpha}^{2\pi - \alpha} 
      G(\tilde \varphi(\theta)) \tilde \varphi'(\theta) \, d\theta \\
  &= -\int_{\alpha}^{2\pi - \alpha} 
      G(\varphi(\theta))\varphi'(\theta) \, d\theta \\
  &= - \frac{|\kappa_0|}4 \sum_{k=0}^3 C_k \int_{\alpha}^{2\pi - \alpha} 
      c^k s^{1-k} \, d\theta
\end{align*}
with the same constants $C_k$. Then, the result \eqref{Ie2:expanded} of the computation below \eqref{Ie2:comp2} yields
\begin{equation*}
  F_\e^2 
  = - \kappa_0 \Big( A_{\phi} \log \frac1{\varepsilon |\kappa_0|} + B_{\phi} + O (\e^2) \Big).   
\end{equation*} 
This completes the proof of \eqref{Fe:t}. 

\begin{rem}[A more precise characterization of $\Psi(x)$] \label{r:Psi}
The proof above shows that
\begin{equation} \label{Psi} 
  \Psi(x)
  = \oint_{\Gamma_\e \cup \cC_\e(x) \cup \gamma_\e(x)} G(y-x) \cdot \tau(y) \, dy + O(\e),    
\end{equation}
where the curves $\cC_\e(x)$ and $\gamma_\e(x)$, the direction in which they are traversed, and the existence of the limit as $\e \to 0$ are detailed in the proof. The proof also shows that the contribution of $\gamma_\e(x)$ to this line integral is $O(\e)$ (uniformly in $x$), and thus its contribution may be left out in \eqref{Psi}. 
\end{rem}

\subsection{Expansion of $\cF_\e$}
\label{s:t:ana:2}

In this section we complete the proof of Theorem \ref{t:ana} by proving \eqref{cFe:t}. The preparation is analogous to that in Section \ref{s:t:ana:1}; we fix an $x \in \Gamma$, translate the coordinate system such that $x$ is at the origin, and remove $x$ from the notation. Again, we set $\kappa_0$, $\tau_0$ and $n_0$ as respectively the curvature, tangent vector (counter-clockwise direction) and outward pointing normal vector of $\Gamma$ at $0$. In addition, in view of the definition of $G_\e$ in \eqref{Ge}, we set
\[
  R_\e(y) := \sqrt{ |y|^2 + \e^2 }
\]
for $y \in \R^2$, and extend this definition to scalars $t \in \R$ by $R_\e(t) = \sqrt{ t^2 + \e^2 }$.
Depending on the sign of $\kappa_0$, we split three cases.

\paragraph{The case $\kappa_0 > 0$.}
Similar to \eqref{Ie123}, we split
\begin{equation}\label{cFe12}
\cF_\e 
= \oint_{\Gamma \cup \cC} G_\e \cdot \tau - \oint_\cC G_\e \cdot \tau
=: \cF_\e^1 - \cF_\e^2,
\end{equation} 
where $\cC$ is traversed in counter-clockwise direction. Since $G_\e$ is regular at $0$, there is no need to avoid the origin, and thus the splitting above is simpler than that in \eqref{Ie123}. Again, we first assume that $\Gamma \cap \cC = \{0\}$, and treat the general case afterwards.

We start with showing that $\cF_\e^1 = \Psi + O(\e)$. We recall $\cC_\e$ and $\gamma_\e$ from Section \ref{s:t:ana:1} (see Figure \ref{fig:Gammae:Ce}), and extend the definition in \eqref{omde} to
\begin{equation*}
  \omega_{\e} 
  := \omega_{0, \e}
  := B_\e(0) \setminus \overline{ B_{r_0}( r_0 n_0) \cup \Omega },
\end{equation*}
i.e., $\omega_\e$ is the union of the two narrow wedges inside $B_\e(0)$ between $\Gamma$ and $\cC$. Then,
\begin{equation} \label{cFe1}
  \cF_\e^1 
  = \oint_{\partial \omega_\e} G_\e \cdot \tau
    + \oint_{\Gamma_\e \cup \cC_\e \cup \gamma_\e} (G_\e - G) \cdot \tau
    + \oint_{\Gamma_\e \cup \cC_\e \cup \gamma_\e} G \cdot \tau.
\end{equation} 
By \eqref{Psi}, the third term equals $\Psi + O(\e)$. Hence, it remains to show that the first two terms are error terms of size $O(\e)$.

We start with the first error term in \eqref{cFe1}. In preparation for applying Stokes' Theorem, we compute
\begin{align} \label{curl:Ge}
  \curl G_\e(y) 
  = \frac{ (1+\nu) y_1^2 + (1 - 2\nu) y_2^2 - (2-\nu) \e^2 }{R_\e^5(y)}
    + \frac32 (1 - \nu) \frac{ \e^2 ( y_1^2 - 4y_2^2 + \e^2) }{R_\e^7(y)}.
\end{align}
Similar to \eqref{g}, it is easy to see that there exists a constant $C > 0$ which only depends on $\nu$ such that 
\[
  \big| \curl G_\e(y) \big| 
  \leq C / |y|^3
  \qquad \text{for all } \e > 0, \: y \in \R^2.
\]
Then, applying Stokes' Theorem to the first error term in \eqref{cFe1}, and then estimating the result analogously to \eqref{OmOme:intd}, we obtain that this error term is $O(\e)$.

To bound the second error term in \eqref{cFe1}, we start with some preparations. We write
\[
  (G_\e - G)(y)^T 
  = \Big( \frac1{ R_\e(y)^3 } - \frac1{ |y|^3 } \Big) y^T \begin{bmatrix}
  0 & -1 \\
  1-\nu & 0
  \end{bmatrix}
  + \frac{3 \varepsilon^2 y^T}{ 2 R_\e(y)^5 } \begin{bmatrix}
  0 & 0 \\
  1 - \nu & 0
  \end{bmatrix}
\]
and claim that $(G_\e - G)(y)$ is small when $y$ remains a fixed distance $\rho > 0$ away from $0$. Indeed, for such $y$,
\[
  \frac1{R_\e(y)} 
  = \frac1{|y|} \frac1{\sqrt{ 1 + \e^2/|y|^2 }}
  = \frac1{|y|} \big( 1 + O(\e^2) \big),
\]
where the term $O(\e^2)$ is uniform in $y$ on $B_\rho(0)^c$, but may depend on $\rho$. Hence, there exists $C > 0$ such that
\[
  \big| (G_\e - G)(y) \big|
  \leq C \e^2 / |y|^2
  \qquad \text{for all } |y| \geq \rho.
\]

To use this result, we expand 
\begin{align} \label{pf:1}
  \oint_{\Gamma_\e \cup \cC_\e \cup \gamma_\e} (G_\e - G) \cdot \tau
  = \bigg( \int_{\Gamma_\rho \cup \gamma_\rho} + \int_{\cC_\rho} + \oint_{\partial \omega_{\e, \rho}} \bigg) (G_\e - G) \cdot \tau,
\end{align}
where $\omega_{\e, \rho}$ is defined in \eqref{omde}.
We claim that the first two integrals are $O(\e^2)$. Indeed, since $\Gamma_\rho \cup \gamma_\rho$ is of finite length, this follows immediately for the first integral. For the second integral, we need to be more precise, because $|\cC_\rho| = O(r_0)$ is not bounded uniformly in $\kappa_0$ (it blows up as $\kappa_0 \downarrow 0$). We use the parametrization of $\cC_\rho$ defined in \eqref{Ie2:vphi} (see \eqref{Ie2:vphi:props} for its properties) to find
\begin{multline*}
  \bigg| \int_{\cC_\rho} (G_\e - G) \cdot \tau \bigg|
  \leq \int_{\rho/r_0}^{2\pi - \rho/r_0} \big| (G_\e - G)(\varphi(\theta)) \big| | \varphi'(\theta) | \, d\theta \\
  \leq C \e^2 \int_{\rho/r_0}^{2\pi - \rho/r_0}  \frac{ | \varphi'(\theta) | }{| \varphi(\theta)) |^2} \, d\theta 
  \leq C \frac{ \e^2 }{4 r_0} \int_{\rho/r_0}^{\pi} \frac{d\theta}{ \sin^2(\theta/2) }  
  \leq C' \frac{ \e^2 }{r_0} \int_{\rho/r_0}^{\pi} \frac{d\theta}{ \theta^2 }
  \leq C'' \e^2, 
\end{multline*}
which is uniform in $\kappa_0$.

Finally, we bound the third integral in \eqref{pf:1}. To apply a similar estimate as for the expansion of $F_\e$, we choose $\rho > 0$ small enough with respect to $\max_{\Gamma} |\kappa|$ such that $\omega_{\e, \rho}$ can be described as in \eqref{OmOme:pmzd}. Then, in preparation for applying Stokes' Theorem, we note that the curl of the integrand is of the form (recalling \eqref{g} and \eqref{curl:Ge})
\[
  \curl (G_\e - G)(y) 
  = \big( C_1 y_1^2 + C_2 y_2^2  \big) \Big( \frac1{R_\e^5(y)} - \frac1{|y|^5} \Big) + \frac{\e^2 (C_3 y_1^2 + C_4 y_2^2 + C_5 \e^2) }{R_\e^7(y)},
\]
where $C_i \in \R$ are constants which depend only on $\nu$. Writing $r := |y|$ and using $R_\e(y) \geq r$, we estimate
\begin{multline*} 
  r^3 \big| \curl (G_\e - G)(y) \big|
  \leq C \Big( 1 - \frac{r^6}{R_\e^6(r)} \Big) + C' \frac{ \e^2 }{R_\e^2(r)} \\
  = C \frac{ 3 (r/\e)^4 + 3 (r/\e)^2 + 1 }{(1 + (r/\e)^2)^3} + \frac{ C' }{1 + (r/\e)^2} 
  =: \psi(r/\e).
\end{multline*}
Note that $\psi$ is integrable on $(0,\infty)$. Then, applying Stokes' Theorem on the third integral in \eqref{pf:1} and recalling \eqref{OmOme:intd}, we obtain
\begin{multline*} 
  \bigg| \oint_{\partial \omega_{\e, \rho}} (G_\e - G) \cdot \tau \bigg|
  = \bigg| \iint_{\omega_{\e, \rho}} \curl (G_\e - G) \bigg| \\
  \leq C \int_\e^\rho r^3 \big| \curl (G_\e - G)(r) \big| \, dr
  \leq C \e \int_1^{\rho/\e} \psi(t) \, dt
  = O(\e).
\end{multline*}
This completes the expansion of the right-hand side of \eqref{cFe1}. Collecting all computations above, we obtain
\[
  \cF_\e^1 = \Psi + O(\e).
\]

Next we expand the second term $\cF_\e^2$ of $\cF_\e$ in \eqref{cFe12}. Inserting \eqref{Ge}, $\cF_\e^2$ reads as
\[
  \cF_\e^2
  = \oint_\cC \frac{y}{ R_\e^3(|y|) } \cdot \begin{bmatrix}
  0 & -1 \\
  1-\nu & 0
  \end{bmatrix} \tau (y) + \frac{3 \varepsilon^2 y}{ 2 R_\e^5(|y|) } \begin{bmatrix}
  0 & 0 \\
  1 - \nu & 0
  \end{bmatrix} \tau (y) \, dy.
\]
Using the parametrization of $\cC$ defined in \eqref{Ie2:vphi} and applying the same steps as in the derivation leading to \eqref{I2e:1}, we obtain
\begin{align*} 
  \cF_\e^2
  = \kappa_0 \Big( 2 \nu (1 - 2 c_\phi^2) I_3^4  
   + (3 \nu c_\phi^2 - \nu - 1) I_3^2  
   + 3 (1 - \nu) (2 c_\phi^2 - 1) I_5^4 
   + \frac32 (1 - \nu) (1 - 3 c_\phi^2) I_5^2 \Big),
\end{align*}
where
\begin{align*}
  I_j^i := \int_{0}^{2\pi} 
      \frac{ 2 r_0^3 \e^{j-3} \sin^i \tfrac\theta2 }{ \sqrt{ 4 r_0^2 \sin^2 \tfrac\theta2 + \e^2  }^j } \, d\theta \qquad \text{for } i = 2,4 \text{ and } j = 3,5.
\end{align*}

To expand the integrals $I_j^i$, we change variables $\theta/2 \to \theta$, use the symmetry of $\sin$ to cut the integration domain in half, and introduce 
\[
  \ell := \frac\e{2 r_0}
\]
to simplify it to
\begin{align*} 
  I_j^i = \ell^{j-3} \int_{0}^{\frac\pi2} 
      \frac{ \sin^i \theta }{ \sqrt{ \sin^2 \theta + \ell^2  }^j } \, d\theta
\end{align*}
By expanding in the integrands the numerator as a polynomial of $\sin^2 \theta + \ell^2$, we obtain
\begin{align*}
  I_3^4
  &= \int_{0}^{\frac\pi2} \sqrt{ \sin^2 \theta + \ell^2  } \, d\theta
     - 2 \ell^2 \int_{0}^{\frac\pi2} \frac1{ \sqrt{ \sin^2 \theta + \ell^2  } } \, d\theta
     + \ell^4 \int_{0}^{\frac\pi2} \frac1{ \sqrt{ \sin^2 \theta + \ell^2  }^3 } \, d\theta \\
  I_5^4
  &= \ell^2 \int_{0}^{\frac\pi2} \frac1{ \sqrt{ \sin^2 \theta + \ell^2  } } \, d\theta
     - 2 \ell^4 \int_{0}^{\frac\pi2} \frac1{ \sqrt{ \sin^2 \theta + \ell^2  }^3 } \, d\theta
     + \ell^6 \int_{0}^{\frac\pi2} \frac1{ \sqrt{ \sin^2 \theta + \ell^2  }^5 } \, d\theta \\   
  I_3^2
  &= \int_{0}^{\frac\pi2} \frac1{ \sqrt{ \sin^2 \theta + \ell^2  } } \, d\theta
     - \ell^2 \int_{0}^{\frac\pi2} \frac1{ \sqrt{ \sin^2 \theta + \ell^2  }^3 } \, d\theta \\
  I_5^2
  &= \ell^2 \int_{0}^{\frac\pi2} \frac1{ \sqrt{ \sin^2 \theta + \ell^2  }^3 } \, d\theta
     - \ell^4 \int_{0}^{\frac\pi2} \frac1{ \sqrt{ \sin^2 \theta + \ell^2  }^5 } \, d\theta. 
\end{align*}
The integrals appearing are complete elliptic integrals of a certain order. To see this, set
\[
  k^2 := \frac1{1 + \ell^2} = \frac{4 r_0^2}{ 4 r_0^2 + \e^2 } \xto{\e \to 0} 1 
\]
and note that
\[
  \sin^2 \theta + \ell^2 = 1 - \cos^2 \theta + \ell^2 = \frac1{k^2} (1 - k^2 \cos^2 \theta).
\]
Hence,
\[
  \int_{0}^{\frac\pi2} \frac1{\sqrt{ \sin^2 \theta + \ell^2 }^m} \, d\theta
  = k^{m} \int_{0}^{\frac\pi2} \frac1{\sqrt{ 1 - k^2 \cos^2 \theta }^m} \, d\theta
\]
for $m \in \Z$, which are complete elliptic integrals when $m$ is a positive, odd integer. For such $m$, the number $m$ is the order of the complete elliptic integral.

Next we expand the appearing complete elliptic integrals around $k=1$. For $m=1$, we obtain from \cite[(8.113.3)]{GradshteynRyzhik14} that
\begin{equation*} 
    \int_0^{\frac\pi2} \frac{d \theta}{ \sqrt{ 1 - k^2 \cos^2 \theta } } 
    = \log \frac4{\sqrt{1 - k^2}} 
      + O \Big( (1 - k^2) \log \frac1{1 - k^2} \Big)
    = \log \frac1\ell + \log 2 + O(\ell^2 |\log \ell|),
\end{equation*} 
and for $m=-1$, \cite[(8.114.3)]{GradshteynRyzhik14} states
\begin{equation*} 
    \int_0^{\frac\pi2} \sqrt{ 1 - k^2 \cos^2 \theta } \, d \theta
    = 1 + O \Big( (1 - k^2) \log \frac1{1 - k^2} \Big)
    = 1 + O (\ell^2 |\log \ell|).
\end{equation*}  
For $m=3$ and $m=5$, we did not find an explicit expansion in the literature. To obtain such expansions, one can either rewrite the integrals in terms of more standard integrals (see e.g.\ \cite[Chap.\ 17]{AbramowitzStegun64}), or replace the cosine by its Taylor polynomial, show that the error made is $O(1)$, and use \cite[(2.271.4--6)]{GradshteynRyzhik14} to expand the obtained integrals. This yields
\begin{align*} 
  \int_{0}^{\frac\pi2} \frac1{\sqrt{ 1 - k^2 \cos^2 \theta }^3} \, d\theta
  &= \frac 1{1 - k^2} + O \Big( \log\frac1{1-k^2} \Big)
  = \frac1{\ell^2} + O (|\log \ell|) \\ 
  \int_{0}^{\frac\pi2} \frac1{\sqrt{ 1 - k^2 \cos^2 \theta }^5} \, d\theta
  &= \frac23 \frac{1}{(1 - k^2)^2} + O \Big( \frac{1}{1 - k^2} \Big)
  = \frac23 \frac1{\ell^4} + O (\ell^{-2}).
\end{align*}

Finally, tracing back our expansions up to $\cF_\e^2$, we observe that
\begin{align*}
  I_3^4 &= 1 + O(\ell^2 |\log \ell|) \\
  I_5^4 &= O(\ell^2 |\log \ell|) \\
  I_3^2 &= \log \frac1\ell + \log 2 - 1 + O(\ell^2 |\log \ell|) \\
  I_5^2 &= \frac13 + O(\ell^2 |\log \ell|) 
\end{align*}
and obtain 
\begin{equation*} 
  \cF_\e^2
  = - \kappa_0 \Big( A_\phi \log \frac1 {\varepsilon \kappa_0} + B_{\phi} + C_{\phi} \Big) + O (\varepsilon),
\end{equation*}
where $A_{\phi}$, $B_{\phi}$ and $C_{\phi}$ are the constants defined in \eqref{ABC}. This completes the proof of Theorem \ref{t:ana} for the case $\kappa_0 > 0$ and $\Gamma \cap \cC = \{0\}$.

Next we demonstrate how this result can be modified to the case of general curves $\Gamma$. The only required modifications are in the definitions of $\omega_\e$ and $\omega_{\e,\rho}$. By taking first $\rho$ and then $\e$ smaller if necessary, we can describe these sets analogously to the description at the end of the case $\kappa_0 > 0$ of the proof for the expansion of $F_\e$. Similar to that proof, we can then apply Stokes' Theorem on $\omega_\e$ and $\omega_{\e,\rho}$ to justify the steps above. This completes the proof of Theorem \ref{t:ana} for the case $\kappa_0 > 0$.

\paragraph{The case $\kappa_0 = 0$.}
As for the expansion of $F_\e$, we can treat the case $\kappa_0 = 0$ similarly to the case $\kappa_0 > 0$ with minor modifications. In fact, two out of three modifications are the same: 
\begin{enumerate}
  \item the replacement of $\cC$ by the tangent line to $\Gamma$ at $0$, and
  \item the observation that $\cF_\e^2 = 0$ (because $G_\e$ is odd).
\end{enumerate}
The third modification is the treatment of the second integral (the one over $\cC_\rho$) in \eqref{pf:1}. By the symmetry of $\cC_\rho$ and the oddness of both $G_\e$ and $G$, we directly obtain 
\[
  \int_{\cC_\rho} (G_\e - G) \cdot \tau = 0.
\]

\paragraph{The case $\tilde \kappa_0 < 0$.} This case can be treated along the same lines as for the expansion of $F_\e$; we omit the details. This completes the proof of Theorem \ref{t:ana}.

\section{Proof of Theorem \ref{t:num}}
\label{s:t:num}

We fix some $i \in \{1,\dots,N\}$, and translate the points $x_j$ of $\Gamma^h$ by $-x_i$. Then, we relabel the points by setting 
\[
  y_j := x_{i+j} - x_i \qquad \text{for all } j \in \Z
\] 
Note that $y_0 = 0$. As in the proof of Theorem \ref{t:ana}, we also translate $\Gamma$ by  $-x_i$, and set $\tau_0$, $n_0$ and $\kappa_0$ as respectively the tangent vector, the normal vector and curvature of the translated copy of $\Gamma$ at $0$. Furthermore, we remove $x$ from the notation, and note that \eqref{Fe:t} reads as
\begin{equation} \label{Fe:t:trld}
  F_\e
  = \kappa_0 \Big( A_{\phi} \log \frac1 {\varepsilon |\kappa_0|} + B_{\phi} \Big) 
  + \Psi + O (\varepsilon).
\end{equation}

We give a constructive proof of \eqref{err:t}, starting with the first of the two terms in the left-hand side. This means that we are going to approximate $F_\e$ as a function of $\by$ without using any further information of $\Gamma$. We refer to this process as a discretization of $F_\e$. As a result, the discretization in \eqref{Feih} will appear as the leading order term in this approximation. 

We start by discretizing the local part (i.e.\ the first term in \eqref{Fe:t:trld}) of $F_\e$. With this aim, it is enough to discretize $\kappa_0$. For later use, we also discretize $\tau_0$, $n_0$ and the constants $A_{\phi}$, $B_{\phi}$ and $C_{\phi}$. We use a simple and standard approximation. Since there are many different discretizations available in the literature on parametric curves, we derive our approximation in detail. 

For this discretization we only use the two points $y_1$ and $y_{-1}$, and denote them in consistency with \eqref{ypm} by respectively $y_+$ and $y_-$. We start with some preliminaries. Let $\varphi$ be the arc length parametrization of $\Gamma$ around $0$ with $\varphi(0) = 0$ and $\varphi'(0) = \tau_0$. Let $t_- < 0 < t_+$ be such that $\varphi(t_\pm) = y_\pm$, and note that 
\begin{equation} \label{tpm:lb}
  |y_\pm| \leq t_\pm. 
\end{equation}
We take $h$ small enough such that the part of $\Gamma$ from $0$ to $y_+$ can be described as the graph of a height function $H$ with respect to the line segment $\gamma_1$, i.e.\ as the graph of
$$
  \eta \mapsto \eta \frac{y_+}{|y_+|} + H(\eta) Q \frac{ y_+}{|y_+|},
  \qquad \text{where } Q := \begin{bmatrix} 0 & 1 \\ -1 & 0 \end{bmatrix}
$$
is the rotation matrix by $90$ degrees in clockwise direction. Since $H(0) = H(|y_+|) = 0$, it follows from the Mean Value Theorem that $H'(\eta_*) = 0$ for some $\eta_* \in (0, |y_+|)$. Since $H \in C^3([0, |y_+|])$, we then have that $|H'(\eta)| \leq C |y_+| \leq C' h$. Hence,
\begin{equation} \label{tp:ub}
 t_+ = \int_0^{|y_+|} \sqrt{1 + H'(\eta)^2 } d\eta 
  \leq |y_+| (1 + C h^2).
\end{equation}
Moreover, since
\[
  H(\eta) 
  = H(0) + \eta H'(0) + O(\eta^2)
  = \eta \big( H'(\eta_*) + O(\eta_*) \big) + O(\eta^2)
  = O(h^2),
\]
we obtain for the region $\omega_1^h$ enclosed by $\Gamma$ and $\gamma_1$ that
\begin{equation} \label{om1h:vol}
  |\omega_1^h| 
  = \int_0^{|y_+|} |H(\eta)| \, d\eta
  \leq C h^3.
\end{equation}
Turning back to \eqref{tp:ub}, one can derive a similar estimate for $t_-$. Together with the lower bound in \eqref{tpm:lb}, this yields
\[
  |y_\pm| \leq |t_\pm| \leq |y_\pm| (1 + C h^2).
\]
We use this to expand $\varphi$ around $0$. Since $\varphi$ is an arc length parametrization, we obtain (see e.g.\ \cite{DoCarmo16}) $\varphi'(0) = \tau_0$ and $\varphi''(0) = \kappa_0 n_0$. Then,
\begin{align} \notag
  y_\pm 
  = \varphi(t_\pm) 
  &= \varphi(0) + t_\pm \varphi'(0)  + \frac12 t_\pm^2 \varphi''(0) + O(t_\pm^3) \\\notag
  &= t_\pm \tau_0  + \frac{1}2 t_\pm^2 \kappa_0 n_0 + O(h^3) \\\label{vphi:Tay}
  &= \pm |y_\pm| \tau_0  + \frac{1}2 |y_\pm|^2 \kappa_0 n_0 + O(h^3).
\end{align}
This completes the preliminaries.

We use \eqref{vphi:Tay} to construct approximations for $\tau_0$, $n_0$ and $\kappa_0$. For any linear combination, we obtain
\begin{equation} \label{ypm:lin:comb}
  a y_+ + b y_- 
  = (a |y_+| - b |y_-|) \tau_0  + \frac{1}2 (a |y_+|^2 + b |y_-|^2) \kappa_0 n_0 + O \big((|a| + |b|)h^3 \big).
\end{equation}
Solving for $a,b$ such that the prefactors of $\tau_0$ and $\kappa_0 n_0$ are respectively $1$ and $0$, we obtain
\begin{align*}   
  a &= \frac{|y_-| / |y_+|}{|y_+| + |y_-|} = O(h^{-1}), \\
  b &= -\frac{|y_+| / |y_-|}{|y_+| + |y_-|} = O(h^{-1}), \\
  \tilde \tau_0^h &:= a y_+ + b y_- = \tau_0 + O(h^2).
\end{align*}
In particular, $|\tilde \tau_0^h| = 1 + O(h^2)$. We use this to normalize our approximation of $\tau_0$:
\begin{equation*} 
  \tau_0^h := \frac{\tilde \tau_0^h}{|\tilde \tau_0^h|} = \tau_0 + O(h^2).
\end{equation*}
Then, we simply rotate to obtain
\begin{equation*}
  n_0^h := Q \tau_0^h = n_0 + O(h^2).
\end{equation*}
Note that this definition of $n_0^h$ is consistent with that in \eqref{nhxi}.
Analogously to $\phi$, we take $\phi^h \in [0, 2\pi)$ such that
\[
  n_0^h = \begin{bmatrix}
    \cos \phi^h \\ \sin \phi^h
  \end{bmatrix}.
\]
Note that with the two components of $n_0^h$, we can approximate the constants $A_{\phi}$, $B_{\phi}$ and $C_{\phi}$ in \eqref{ABC} by respectively $A_{\phi^h}$, $B_{\phi^h}$ and $C_{\phi^h}$ with an error of size $O(h^2)$.

Similarly, solving for $a,b$ such that the prefactors of $\tau_0$ and $\kappa_0 n_0$ in \eqref{ypm:lin:comb} are respectively $0$ and $1$, we obtain
\begin{equation*} 
  a = \frac{2 / |y_+|}{|y_+| + |y_-|} = O(h^{-2}), \qquad
  b = \frac{2 / |y_-|}{|y_+| + |y_-|} = O(h^{-2}), \qquad
  a y_+ + b y_- = \kappa_0 n_0 + O(h).
\end{equation*}
Multiplying both sides of the third equation by $\tilde n_0^h = Q \tilde \tau_0^h$, we obtain
\begin{equation} \label{kaph:expa}
  \kappa_0^h := (a y_+ + b y_-) \cdot \tilde n_0^h = \kappa_0 + O(h).
\end{equation}
Rewriting this in terms of $y_\pm$, we get
\begin{equation*} 
  \kappa_0^h 
  = 2 \Big( \frac1{|y_+|} + \frac1{|y_-|} \Big) \frac{ y_- \cdot Q y_+ }{ (|y_+| + |y_-|)^2 },
\end{equation*}
which motivates the expression in \eqref{kaphxi}.
 
Finally, we substitute the expansions above in the local term of \eqref{Fe:t:trld}. While this needs no further motivation for most of the terms, we wish to treat the expansion of $\psi(\kappa_0) := \kappa_0 \log |\kappa_0|$ carefully. From the derivation of \eqref{kaph:expa} we obtain that
\[  
  |\underbrace{ \kappa_0 - \kappa_0^h }_{R^h}| \leq M h
\]
for a constant $M > 0$ which is independent of $h$ and $\kappa_0$. If $|\kappa_0| \leq 2 M h$, then we observe that $\psi(\kappa_0) = O(h |\log h|)$ and $\psi(\kappa_0^h) = O(h |\log h|)$. If $|\kappa_0| \geq 2 M h$, then $M h \leq \kappa_0^h \leq C$, and thus we may apply Taylor's Theorem on $\psi$ at $\kappa_0^h$. This yields an $\kappa^* \in \overline{B(\kappa_0, Mh)}$ such that
\[
  \psi(\kappa_0)
  = \psi(\kappa_0^h) + R^h \psi'(\kappa^*)
  = \psi(\kappa_0^h) + R^h (\log |\kappa^*| + 1)
  = \psi(\kappa_0^h) + O(h |\log h|).
\]
Using this, we obtain from \eqref{Fe:t:trld} that
\begin{equation} \label{Fe:t:trld:1}
  F_\e
  = \kappa_0^h \Big( A_{\phi^h} \log \frac1 {\varepsilon |\kappa_0^h|} + B_{\phi^h} \Big) 
  + \Psi + O \big( \varepsilon + h |\log \e| + h |\log h| \big).  
\end{equation}  

It is left to approximate the nonlocal term $\Psi$. We do this in detail for the case $\kappa_0 > 0$; the case $\kappa_0 \leq 0$ can then be treated along the same lines as in the proof of Theorem \ref{t:ana}. 

We start with some preliminaries. Let $\varphi$ be again the arc length parametrization of $\Gamma$, and let $t_j$ be defined by
\[ \varphi(t_j) = y_j \qquad \text{for } - \lfloor  N/2 \rfloor \leq j \leq \lfloor  N/2 \rfloor - 1.  \] 
Let $\Gamma_j$ be the part of $\Gamma$ in between $y_{j-1}$ and $y_j$, and set $\omega_j^h$ as the region enclosed by $\Gamma_j$ and $\gamma_j$. Analogously to \eqref{om1h:vol}, one can derive (for $h$ small enough) that 
\begin{equation} \label{omjh:vol}
  |\omega_j^h| \leq C h^3 \qquad \text{for all } j \in \Z.
\end{equation}

Interpreting $d(s) := |\varphi(s)|$ as the distance from the origin, we note that
\[
  d'(s) = \frac{\varphi(s)}{|\varphi(s)|} \cdot \varphi'(s)
  \qquad \text{for all } s \in \R \setminus |\Gamma| \Z.
\]
Since $\varphi(s) / |\varphi(s)| \to \tau_0$ as $s \downarrow 0$, $d'(s)$ is uniformly continuous on $(0, \rho]$ for any $\rho < |\Gamma|$. Observing that $d'(0\pm) = \pm1$, we may take $\rho$ independent of $h$ such that
\begin{equation} \label{dp:ests}
\inf_{(0, \rho)} d' \geq \frac12
   \qand
   \sup_{(-\rho,0)} d' \leq - \frac12,
\end{equation}
and such that $\Gamma$ intersects $\partial B_\rho(0)$ in precisely two points. Let $n^\rho$ be the largest integer for which
\[
  |y_{n^\rho}| \leq \rho 
  \qand 
  |y_{-n^\rho}| \leq \rho.
\]
We observe from
\[
  |y_{n^\rho}| 
  = |y_{n^\rho} - y_0| 
  \leq \sum_{j=1}^{n^\rho} |y_j - y_{j-1}|
  \leq C n^\rho h
\]
that $n^\rho \geq c / h$ for some $c > 0$ independent of $h$. Recalling $m^h$ from \eqref{mh}, we take $h$ small enough such that $m^h \leq n^\rho$.

From this construction, we obtain
\begin{multline} \label{yj:norm:incr}
  |y_{j}| - |y_{j-1}| 
  = d(t_j) - d(t_{j-1})
  = \int_{t_{j-1}}^{t_j} d'(s) \, ds \\
  \geq \frac12 (t_j - t_{j-1})
  \geq \frac12 |\gamma_j|
  \geq c h
  \qquad \text{for } j = 1,\ldots,n^\rho.
\end{multline}
Since a similar estimate holds for the points $y_{-j}$, we obtain
\begin{align} \label{ypm:mh:est}
  c h^{2/3} \leq |y_{\pm m^h}| \leq C h^{2/3}
\end{align} 
for all $h$ small enough.
For convenience, we first assume that 
\begin{align}  \label{eh:assns}
  \e^h := |y_{-m^h}| = |y_{m^h}|,
\end{align}
and comment on the general case afterwards. Note from \eqref{ypm:mh:est} that
\[
  \e^h = O(h^{2/3}).
\]
This concludes the preliminaries for discretizing $\Psi$ in \eqref{Fe:t:trld:1}.

From the characterization of $\Psi$ in Remark \ref{r:Psi}, we obtain
\begin{align} \label{Stokes:h}
\Psi
= \int_{\Gamma_{\e^h}} G \cdot \tau + \int_{\cC_{\e^h}} G \cdot \tau + O(\e_h)
\end{align}  
for any $h > 0$ small enough with respect to $\Gamma$. We discretize the first term by the integral over 
\begin{equation} \label{Gammaehh}
 \Gamma_{\e^h}^h 
  := \Gamma^h \setminus B(0, \e^h) 
  = \bigcup_{j = m^h + 1}^{N - m^h} \gamma_j 
  \subset \Gamma^h.
\end{equation}
Figure \ref{fig:Gammaeh} illustrates the setting.

\begin{figure}[ht]
\centering
\begin{tikzpicture}[scale=1.5, >= latex]
\def \rr {.04} 

\begin{scope}[xscale = 2]  
  \begin{scope}[rotate = 240]
    \draw[dotted] (0,0) circle (1);
    \draw[thick] (1,0) arc (0:240:1);
  \end{scope}  
\end{scope}

\fill (2,0) circle (\rr);
\fill (0,1) circle (\rr);
\fill (0,-1) circle (\rr);

\foreach \x in {0,180}{ 
\foreach \pmm in {1,-1}{ 
  \begin{scope}[rotate = \x]  
    \fill (\pmm*1.732,.5) circle (\rr);
    \fill (\pmm*1,.866) circle (\rr);  
    \draw[red, thick] (\pmm*2,0) -- (\pmm*1.732,.5) -- (\pmm*1,.866) -- (0,1);   
  \end{scope}
}
}

\draw[white, very thick] (-1,-.866) -- (-1.732,-.5) -- (-2,0) -- (-1.732,.5) -- (-1,.866);
\draw[dotted, red] (-1,-.866) -- (-1.732,-.5) -- (-2,0) -- (-1.732,.5) -- (-1,.866); 
\fill[black!15!white] (-1.732,.5) circle (\rr);
\fill[black!15!white] (-1.732,-.5) circle (\rr);
\fill[black] (-2,0) circle (\rr);
\fill[black] (-1,.866) circle (\rr);
\fill[black] (-1,-.866) circle (\rr);

\draw[green!70!black, dotted] (-2,0) circle (1.323);
\draw[green!70!black] (-3.323,0) node[right]{$B_{\e^h}(0)$};

\draw (-2,0) node[left] {$0$};
\draw (1.732,-.5) node[anchor = north west] {$y_{j-1}$};
\draw (2,0) node[right] {$y_j$};
\draw (1.9,-.35) node[right] {$\Gamma_j$};
\draw (-1,.866) node[above] {$y_{-m^h}$};
\draw (-1,-.866) node[below] {$y_{m^h}$};
\draw[red] (1.732,-.15) node {$\gamma_j$};
\draw[red] (.5,.93) node[below] {$\Gamma_{\e^h}^h$};
\draw (.5,1) node[above] {$\Gamma_{\e^h}$};
\end{tikzpicture} \\
\caption{$\Gamma_{\e^h}$ and $\Gamma_{\e^h}^h$ related to $\Gamma$ and $\Gamma^h$ as illustrated in Figure \ref{fig:Gammah}. The sketch is a special case in which (\ref{eh:assns}) is satisfied.}
\label{fig:Gammaeh}
\end{figure}

Then, by Stokes' Theorem\footnote{In the case where $\gamma_j$ and $\Gamma_j$ intersect, we apply a similar argument as that at the end of the case $\kappa_0 > 0$ in the proof of \eqref{Fe:t}.}
\begin{equation} \label{pf:3}
\bigg| \int_{\Gamma_{\e^h}} G \cdot \tau - \int_{\Gamma_{\e^h}^h} G \cdot \tau \bigg|
\leq \sum_{j = m^h + 1}^{N-m^h} \bigg| \oint_{\gamma_j \cup \Gamma_j} G \cdot \tau \bigg|
\leq \sum_{j = m^h + 1}^{N-m^h} \iint_{\omega_j^h} |g|,
\end{equation}
where $g = \curl G$ (see \eqref{g}).
To bound the sum in the right-hand side, we split it in three sub-summations. For $j = n^\rho + 1, \ldots, N - n^\rho$, by the construction of $\rho$, it holds that $\dist(0, \Omega_j) \geq c$. Hence, 
\[ \max_{\omega_j^h} |g| \leq C. \]
Using this together with \eqref{omjh:vol}, we obtain
\[
  \iint_{\omega_j^h} |g|
  \leq |\omega_j^h| \max_{\omega_j^h} |g| 
  \leq C h^3.
\]
Noting from \eqref{h:xi} that $N \leq C/h$, we then have that
\[
  \sum_{j = n^\rho + 1}^{N-n^\rho} \iint_{\omega_j^h} | g |
  \leq C N h^3
  \leq C' h^2.
\]

The remaining two sub-summations in the right-hand side of \eqref{pf:3} can be treated similarly to one another; we focus on the one over $j = m^h + 1,\ldots ,n^\rho$. From the preliminaries (see \eqref{dp:ests} and \eqref{yj:norm:incr}) we obtain
\[
  \max_{\omega_j^h} |g|
  \leq \max_{y \in \omega_j^h} \frac C{|y|^3}
  = \frac C{|y_{j-1}|^3} 
  \leq \frac{C'}{(hj)^3}.
\]
Then, 
\[
  \sum_{j = m^h + 1}^{n^\rho} \bigg| \iint_{\omega_j^h} g \bigg|
  \leq \sum_{j = m^h + 1}^{n^\rho} |\omega_j^h| \max_{\omega_j^h} |g|
  \leq \sum_{j = m^h + 1}^{n^\rho} \frac{C}{j^3}
  \leq \int_{m^h}^\infty \frac{C}{\alpha^3} \, d\alpha
  = 2C / (m^h)^2
  \leq C' h^{2/3}.
\]
In conclusion, recalling \eqref{pf:3} and \eqref{Gammaehh},
\begin{equation} \label{pf:8}
\int_{\Gamma_{\e^h}} G \cdot \tau 
= \sum_{j = m^h + 1}^{N-m^h} \int_{\gamma_i} G \cdot \tau  + O(h^{2/3}).
\end{equation}

Next, we estimate the second term in \eqref{Stokes:h} by the integral over a circle $\cC^h$ which we will construct from $n_0^h$ and $\kappa_0^h$. For $\cC^h$ to be close to $\cC$, we require in \eqref{kaph:expa} that the error term $O(h)$ is sufficiently smaller than $\kappa_0$. This motivates us to first consider the case $\kappa_0 \geq \e^h$, and treat the case for small $\kappa_0$ afterwards.  

Assuming that $\kappa_0 \geq \e^h$, we set 
\[
  \cC^h := \partial B_{r_0^h} (r_0^h n_0^h) 
  \qand
  \cC_{\e^h}^h := \cC^h \setminus B_{\e^h}(0),
\]
where
\[
  r_0^h = \frac1{\kappa_0^h} = \frac1{\kappa_0 + O(h)} = \frac{1 + O(h/\kappa_0)}{\kappa_0}
  = r_0 (1 + O(r_0 h) )
\] 
can be computed from the two points $y_+$ and $y_-$. 
Let $\gamma_{\e^h}$ be the two small arcs on $\partial B_{\e^h}(0)$ which connect the endpoints of $\cC_{\e^h}^h$ and $\cC_{\e^h}$, and let $\omega^h$ be the region enclosed by the closed loop $\cC_{\e^h}^h \cup \cC_{\e^h} \cup \gamma_{\e^h}$; see Figure \ref{fig:Cehh} for a sketch. Then, by Stokes' Theorem
\begin{equation} \label{pf:4}
  \bigg| \int_{\cC_{\e^h}} G \cdot \tau - \int_{\cC_{\e^h}^h} G \cdot \tau \bigg|
  \leq \bigg| \int_{\gamma_{\e^h}} G \cdot \tau \bigg| + \iint_{\omega^h} |g|.
\end{equation}

\begin{figure}[ht]
\centering
\begin{tikzpicture}[scale=2, >= latex]
\def \rr {.03} 

\begin{scope}[rotate = 195] 
\draw [thick, blue] (1,0) circle (1);
\fill [blue] (1,0) circle (\rr);
\end{scope}
\begin{scope}[scale = .9, rotate = 165] 
\draw [thick, red] (1,0) circle (1);
\fill [red] (1,0) circle (\rr);
\end{scope}

\fill [white] (0,0) circle (.4);
\draw[green!70!black, dotted] (0,0) circle (.4);


\begin{scope}[rotate = 195] 
\draw [dotted, blue] (1,0) circle (1);
\draw[blue, <->] (0,.8) node[right]{$\tau_0$} -- (0,0) -- (.8,0) node[anchor = north west]{$n_0$};
\draw[blue] (2,0) node[left]{$\cC_{\e^h}$};
\begin{scope}[shift={(1,0)}, rotate = 60] 
  \draw[blue] (0,0) -- (1,0) node[midway, left] {$r_0$};    
\end{scope}
\end{scope}
\begin{scope}[rotate = 165] 
\begin{scope}[scale = .9]
\draw [dotted, red] (1,0) circle (1);
\draw[red] (2,0) node[left]{$\cC_{\e^h}^h$};
\begin{scope}[shift={(1,0)}, rotate = 45] 
  \draw[red] (0,0) -- (1,0) node[midway, above] {$r_0^h$};    
\end{scope}
\end{scope}
\draw[red, <->] (0,.8) node[left]{$\tau_0^h$} -- (0,0) -- (.8,0) node[anchor = south west]{$n_0^h$};
\end{scope}

\begin{scope}[scale = .4, rotate = 87] 
  \draw[green!70!black, thick] (1,0) arc (0:28:1);       
\end{scope}
\begin{scope}[scale = .4, rotate = 240] 
  \draw[green!70!black, thick] (1,0) arc (0:35:1);       
\end{scope}

\draw[green!70!black] (.4,0) node[right]{$B_{\e^h}(0)$};
\draw[green!70!black] (.5,.8) node {$\gamma_{\e^h}$};
\begin{scope}[shift={(.5,.8)}, rotate = 213] 
  \draw[->, green!70!black] (0.15,0) -- (.7,0);    
\end{scope}
\begin{scope}[shift={(.5,.8)}, rotate = 245] 
  \draw[->, green!70!black] (0.15,0) -- (1.28,0);    
\end{scope}
\fill (0,0) circle (\rr);
\end{tikzpicture} \\
\caption{Sketch of the closed loop $\cC_{\e^h}^h \cup \cC_{\e^h} \cup \gamma_{\e^h}$.}
\label{fig:Cehh}
\end{figure} 

Next we show that both integrals in the right-hand side of \eqref{pf:4} are small. We start with the first one. We parametrize $\gamma_{\e^h}$ by 
\begin{equation} \label{varphi:s}
  \varphi(\theta) := \e^h \begin{bmatrix} \cos \theta \\ \sin \theta \end{bmatrix} \qquad \text{for } \theta \in \Theta,
\end{equation}
where $\Theta$ is (similar to \eqref{OmOme:pmzd}) the union of two intervals. Let $\theta_1$ and $\theta_2$ be the endpoints of one of these intervals. In particular, $\varphi(\theta_1) \in \cC$ and $\varphi(\theta_2) \in \cC^h$, i.e.
\[
  |\varphi (\theta_1) - r_0 n_0| = r_0 \qand |\varphi (\theta_2) - r_0^h n_0^h| = r_0^h.
\]
To solve for $\theta_1$, we first compute
\begin{multline*}
  r_0^2 
  = |\varphi (\theta_1) - r_0 n_0|^2
  = (\e^h \cos \theta_1 - r_0 \cos \phi)^2 + (\e^h \sin \theta_1 - r_0 \sin \phi)^2 \\
  = (\e^h)^2 + r_0^2 - 2 \e^h r_0 \cos (\theta_1 - \phi).
\end{multline*}
Then, we obtain two solutions given by
\[
  \theta_1 = \phi \pm \arccos \frac{\e^h}{2r_0}.
\]
Analogously, we obtain
\[
  \theta_2 = \phi^h \pm \arccos \frac{\e^h}{2r_0^h}.
\]
Taking the plus sign in both equations above, we obtain the endpoints of one of the two intervals of $\Theta$. Then, substituting the expansions for $\phi^h$ and $\kappa_0^h = 1/r_0^h$, we obtain
\begin{equation} \label{Thetas:size}
  \big| \theta_2 - \theta_1 \big|
  \leq C h^2 + \Big| \arccos \frac{\e^h}{2r_0} - \arccos \Big( \frac{\e^h}{2r_0} + O({\e^h} h) \Big) \Big|
  \leq C' {\e^h} h.
\end{equation}
Hence, $|\Theta| \leq C {\e^h} h$. We use this to estimate the first integral in \eqref{pf:4} by
\[
  \bigg| \int_{\gamma_{\e^h}} G \cdot \tau \bigg|
  = \bigg| \int_{\Theta} G( \varphi(\theta) ) \cdot \varphi'(\theta) \, d\theta \bigg|
  \leq |\Theta| \frac C{\e^h} 
  \leq C' h.
\]

To estimate the second integral in \eqref{pf:4}, we split $\omega^h$ in two pieces; the part inside $B_{r_0}(0)$ and the part outside $B_{r_0}(0)$. For the inside part, we write similar to \eqref{omde}
\begin{align} \label{pf:45}
  \omega^h \cap B_{r_0}(0)
  = \{ (s, \theta) : \e^h < s < r_0, \ \theta \in \Theta(s) \},
\end{align}
where $\Theta(s)$ is the extension of $\Theta$ in \eqref{varphi:s} for $s = \e^h$ to $s \in (\e^h, r_0)$. In particular, the same argument yields $|\Theta(s)| \leq C sh$. Then,
\begin{multline} \label{omplush}
  \iint_{\omega^h \cap B_{r_0}(0)} |g|
  = \int_{\e^h}^{r_0} \int_{\Theta(s)} \big| g(s, \theta) \big| s \, d\theta ds \\
  \leq C \int_{\e^h}^{r_0} |\Theta(s)| \frac1{s^2} \, ds
  \leq C' h \int_{\e^h}^{r_0} \frac1{s} \, ds
  = O(h |\log h|). 
\end{multline}

For the part of $\omega^h$ outside $B_{r_0}(0)$, note that $\omega^h$ remains inside the tubular neighborhood of $\cC = \partial B(r_0 n_0, r_0)$ of size $O(r_0^2 h)$. Indeed, for any point $y \in \cC^h = \partial B(r_0^h n_0^h, r_0^h)$, we obtain from the triangle inequality that
\begin{align*} 
  |y - r_0 n_0| &\leq |y - r_0^h n_0^h| + |r_0^h n_0^h - r_0 n_0| = r_0 + O(r_0^2 h), \\
  |y - r_0 n_0| &\geq |y - r_0^h n_0^h| - |r_0^h n_0^h - r_0 n_0| = r_0 + O(r_0^2 h). 
\end{align*} 
Hence, $|\omega^h| \leq C r_0^3 h$, and thus
\begin{equation} \label{ChCeh:outside:Brho}
  \iint_{\omega^h \setminus B_{r_0}(0)} |g| 
  \leq |\omega^h| \max_{B_{r_0}(0)^c} |g|
  \leq C h.
\end{equation}

Inserting our findings above in \eqref{pf:4}, we obtain
\begin{equation} \label{pf:5}
  \int_{\cC_{\e^h}} G \cdot \tau 
  = \int_{\cC_{\e^h}^h} G \cdot \tau + O(h |\log h|).
\end{equation}
To expand the integral in the right-hand side, we parametrize $\cC_{\e^h}^h$ as in \eqref{Ie2:vphi} by
\begin{equation*}
  \varphi^h (\theta) := r_0^h \begin{bmatrix} \cos \phi^h - \cos (\theta + \phi^h) \\ \sin \phi^h - \sin (\theta + \phi^h) \end{bmatrix} \qquad \text{with } \alpha^h < \theta < 2\pi - \alpha^h,
\end{equation*}
where $\alpha^h \in (0, \pi)$ is such that $|\varphi^h(\alpha^h)| = \e^h$. An analogous derivation as the one for $F_e^2$ leading to \eqref{Ie2:expanded} yields
\begin{equation} \label{pf:55}
  \int_{\cC_{\e^h}^h} G \cdot \tau
  = -\kappa_0^h \Big( A_{\phi^h} \log \frac{1}{\e^h \kappa_0^h} + B_{\phi^h}  \Big) + O \big( (\e^h)^2 \big). 
\end{equation}
Then, \eqref{pf:5} yields 
\begin{align} \label{pf:6}
\int_{\cC_{\e^h}} G \cdot \tau 
  = -\kappa_0^h \Big( A_{\phi^h} \log \frac{1}{\e^h \kappa_0^h} + B_{\phi^h}  \Big) + O(h |\log h|).
\end{align}

We recall that \eqref{pf:6} relies on the assumption $\kappa_0 \geq \e^h$. The remaining case $0 < \kappa_0 < \e^h$ can be treated with minor modifications to the derivation of \eqref{pf:6}. We list these modifications first for the additional assumption $\kappa_0^h > 0$, for which no change to the definition of $\cC^h$ is required. The main modification is that we use $B(0, 1/\e^h)$ instead of $B_{r_0}(0)$ when splitting $\omega_h$ in two pieces. It is easy to see that also $\partial B(0, 1/\e^h)$ intersects with $\cC^h$ in two points, and that the estimate $|\Theta (s)| \leq Csh$ remains valid. Using this estimate, as in \eqref{omplush}, we obtain
\[
  \iint_{\omega^h \cap B(0, 1/\e^h)} |g|
  = O(h |\log h|). 
\]
As an alternative to \eqref{ChCeh:outside:Brho}, we use the rougher estimate
\[
  \iint_{\omega^h \setminus B(0, 1/\e^h)} |g| 
  \leq \iint_{B(0, 1/\e^h)^c} |g|
  \leq C \int_{1/\e^h}^\infty \frac1{s^2} \, ds
  = C \e^h.
\]
Then, the same steps leading to \eqref{pf:6} yield
\begin{align} \label{pf:7}
\int_{\cC_{\e^h}} G \cdot \tau 
  = -\kappa_0^h \Big( A_{\phi^h} \log \frac{1}{\e^h \kappa_0^h} + B_{\phi^h}  \Big) + O(\e^h).
\end{align}

Next we treat the case $0 < \kappa_0 < \e^h$ with $\kappa_0^h < 0$. Note from \eqref{kaph:expa} that this setting implies $\kappa_0 < Ch$ and $\kappa_0^h \geq -C' h$. We set $r_0^h := -1/\kappa_0^h > 0$, $\cC^h := \partial B(- r_0^h n_0^h, r_0^h)$ and
\[
  \omega^h = \big( B(r_0 n_0, r_0) \cup B(- r_0^h n_0^h, r_0^h) \big)^c \cup \big( B(r_0 n_0, r_0) \cap B(- r_0^h n_0^h, r_0^h) \big).
\] 
In addition to the modification in the case $\kappa_0^h > 0$, the derivation of $|\Theta| \leq C \e^h h$ requires a modification too. Indeed, while the endpoint $\theta_1$ can be found analogously, the condition for $\theta_2$ becomes 
\[
  |\varphi (\theta_2) + r_0^h n_0^h| = r_0^h.
\]
Solving for $\theta_2$ yields
\[
  \theta_2 
  = \phi^h \pm \arccos \frac{-\e^h}{2r_0^h}.
\]
Then, using the expansion
\[
  \frac{-1}{r_0^h} = \kappa_0^h = \kappa_0 + O(h),
\] 
we get
\[
  \theta_2 
  = \phi^h \pm \arccos \Big( \frac{\e^h}{2r_0} + O(\e^h h) \Big).
\]
Then, $|\Theta| \leq C \e^h h$ follows from \eqref{Thetas:size} as before.

Finally, we treat the case $0 < \kappa_0 < \e^h$ with $\kappa_0^h = 0$. Setting $\cC^h = \{ t \tau^h : t \in \R \}$, this case can be treated similarly as either the case $\kappa_0^h > 0$ or the case $\kappa_0^h < 0$; we omit the details. This completes the proof of \eqref{pf:7} for all $\kappa_0 > 0$ without any additional assumptions on the sign or size of $\kappa_0^h$. 

In conclusion, starting from \eqref{Fe:t:trld:1} and substituting consecutively \eqref{Stokes:h}, \eqref{pf:8} and \eqref{pf:7}, we obtain
\begin{align} \label{pf:9}
  F_\e 
  = \kappa_0^h A_{\phi^h} \log \frac{\e^h}{\e} 
  + \sum_{j = m^h + 1}^{N-m^h} \int_{\gamma_i} G \cdot \tau + O \big( \varepsilon + h |\log \e| + h^{2/3} \big).
\end{align} 
This prove \eqref{err:t} for the first term in the left-hand side of \eqref{err:t} under the additional assumption $|y_{-m^h}| = |y_{m^h}|$. 

The proof above easily extends to the generic case
\[
  \e^h := |y_{m^h}| \neq |y_{-m^h}| =: \e^{-h}.
\]
Indeed, the main modification is that $B_{\e^h}(0)$ is replaced by the union of two half-balls cut along $n_0$:
\[
  D := \big\{ x \in B_{\e^h}(0) : x \cdot \tau_0 \geq 0 \big\}
  \cup
  \big\{ x \in B_{\e^{-h}}(0) : x \cdot \tau_0 \leq 0 \big\}.
\]
Indeed, in the proof above, the narrow wedges (see, e.g.\ \eqref{pf:45}) between any of the four curves $\Gamma$, $\Gamma^h$, $\cC$ and $\cC^h$ can be treated independently, and are always included in one of the two half-balls. A ramification of this modification is that \eqref{pf:55} changes to 
\begin{align*}
  \int_{\cC^h \setminus D} G \cdot \tau
  &= -\kappa_0^h \Big( \frac12 A_{\phi^h} \log \frac{1}{\e^h \e^{-h} (\kappa_0^h)^2} + B_{\phi^h}  \Big) + O \big( h^{4/3} \big) \\
  &= -\kappa_0^h \Big( \frac12 A_{\phi^h} \log \frac{1}{|y_{m^h}| |y_{-m^h}| (\kappa_0^h)^2} + B_{\phi^h}  \Big) + O \big( h^{4/3} \big);
\end{align*}
this can be seen from obvious modifications to the argument leading to \eqref{Ie2:expanded}. 

To complete the proof of Theorem \ref{t:num}, we note that \eqref{err:t} follows almost directly from \eqref{pf:9}. Indeed, by the triangle inequality and the definitions of $\cF_\e$ and $\cF_{\e,i}^h$, we get
\[
  \big| \cF_{\e,i}^h (\bx) - \cF_\e(x_i) \big| 
  \leq \big| F_{\e,i}^h (\bx) - F_\e(x_i) \big| 
       + \big| \kappa_0^h C_{\phi^h} - \kappa_0 C_\phi \big|.
\]
Recalling from the proof that $\kappa_0^h = \kappa_0 + O(h)$ and $C_{\phi^h} = C_{\phi} + O(h^2)$, we obtain \eqref{err:t}.

\section*{Acknowledgements}

The author gratefully acknowledges support from JSPS KAKENHI Grant Number JP20K14358.
He also expresses his sincere gratitude to Riccardo Scala for the conception of the proof of Theorem \ref{t:ana}. 

\newcommand{\etalchar}[1]{$^{#1}$}

\end{document}